\documentclass[11pt,leqno]{article}
\usepackage{amsthm,amsmath,amsfonts,latexsym,amssymb}
\usepackage{mathrsfs,dutchcal,arydshln,gensymb,appendix,abstract}
\usepackage[dvips]{graphics,color}
\usepackage{graphicx}
\usepackage{tikz}
\date{}

\begin{document}
\title{On Euler-Dierkes-Huisken variational problem}
\author{\centerline {Hongbin Cui and Xiaowei Xu\footnote{The corresponding author}}}
 
\maketitle

\begin{abstract}
In this paper, we study the stability and minimizing properties of higher codimensional surfaces in Euclidean
space associated with the $f$-weighted area-functional 
$$\mathcal{E}_f(M)=\int_M f(x)\; d \mathcal{H}_k$$
with the density function $f(x)=g(|x|)$ and $g(t)$ is non-negative, which develop
the recent works by U. Dierkes and G. Huisken (Math. Ann., 20 October 2023) on hypersurfaces with the density
function $|x|^\alpha$. Under suitable assumptions on $g(t)$, we prove that minimal cones with globally flat normal bundles are $f$-stable, and we also prove that the regular minimal cones satisfying Lawlor curvature criterion, the highly singular determinantal varieties and  Pfaffian varieties without some exceptional cases are $f$-minimizing. As an application, we show that $k$-dimensional minimal cones over product of spheres are $|x|^\alpha$-stable for $\alpha\geq-k+2\sqrt{2(k-1)}$, the oriented stable minimal hypercones are $|x|^\alpha$-stable for 
$\alpha\geq 0$, and we also show that the minimal cones over product of spheres $\mathcal{C}=C \left(S^{k_1} \times \cdots \times S^{k_{m}}\right)$ are $|x|^\alpha$-minimizing for $\dim \mathcal{C} \geq 7$, $k_i>1$ and $\alpha \geq 0$, 
the Simons cones $C(S^{p} \times S^{p})(p\geq 1)$ are $|x|^\alpha$-minimizing for any $\alpha \geq 1$, 
which relaxes the assumption $1\leq\alpha \leq 2p$ in \cite{DH23}.
\end{abstract}

\medskip\noindent


		
\bigskip
		
\tableofcontents
		
\section{Introduction}
\subsection{Backgrounds and related works}

\medskip\noindent

Let $f(x)$ be a function on $\mathbb{R}^{n}$, U. Dierkes and G. Huisken (\cite{DH23}) considered the following 
$f$-\emph{weighted area-functional}
\begin{equation}\label{1-1}
\mathcal{E}_f(M)=\int_M f(x)\; d \mathcal{H}_{n-1}
\end{equation}
for hypersurface $M$ in $\mathbb{R}^{n}$, where $\mathcal{H}_{n-1}$ stands for the $(n-1)$-dimensional Hausdorff measure. They focus on the $f$-stationary, $f$-stable and $f$-minimizing problems for the case that $f(x)=|x|^\alpha$ with $\alpha\in\mathbb{R}$. This variational problem, especially for that $ \alpha>0$, has a long history dating back to the work of L. Euler (\cite{Eu44}) on minimizing the polar moment of inertia, as introduced in \cite{DH23}. They showed that the Euler-Lagrange equation of this functional has three families of canonical solutions: hyperplanes through the origin, hyperspheres, and minimal hypercones. In other words, these solutions are $f$-stationary. U. Dierkes and G. Huisken proved hyperspheres are $f$-stable. They provided the conditions for the $f$-stability of minimal cones and studied the related Bernstein problem. Moreover, they proved that Simons cones $C(S^{\ell} \times S^{\ell})$ are $f$-minimizing cones for $1 \leq \alpha \leq 2\ell$ in the sense of the least $\alpha$-energy, inspired by the work in \cite{BDGG69}. They also proved that hyperplanes and hyperspheres are minimizers of $|x|^{\alpha}$ with different $\alpha$. There are some advancements in finding stationary, stable and minimizing surfaces related to the variational problem of Euler-Dierkes-Huisken type, like $f(x)=x_{n}^{\alpha}$ on $\mathbb{R}^{n-1} \times \mathbb{R}^{+}$, can be found in \cite{D88}, \cite{D89}, \cite{DHT10}, and \cite{Pe20}.
Some other variational problems of hypersurfaces for generalized area-functional, such as the parametric elliptic variational problem, were studied in \cite{AlmSS77},  \cite{M91}, \cite{MY21}, \cite{CL23}, \cite{MY24}.

\medskip\noindent

In this paper, we wish to develop the theory of U. Dierkes and G. Huisken to the higher codimensional case. Namely, we consider the following $f$-weighted area-functional 
\begin{equation}\label{fequation}
 \mathcal{E}_f(M)=\int_M f(x)\; d\mathcal{H}_{k},   
\end{equation}
where $f(x)$ is a continuous function on $\mathbb{R}^n$ and $M\subset\mathbb{R}^n$ is a $k$-dimensional surface.
Here $f(x)$ is also called \emph{density function} in \cite{M05}, \cite{M16}. The functional (\ref{fequation}) will be cited as
\emph{Euler-Dierkes-Huisken functional} and the corresponding variational problem is called \emph{Euler-Dierkes-Huisken variational problem}. Throughout this paper, we always consider a more general density function $f(x)$ instead of $|x|^\alpha$ in \cite{DH23}. That is, we take $f(x)=g(|x|)$, which is a radial function with $g(t)$ is nonnegative and non-decreasing on $\mathbb{R}^+$. If $f(x)$ is $C^1$ on $\mathbb{R}^n$, the Euler-Lagrange equation of \eqref{fequation} is given by 
\begin{equation}\label{1-3}
f \textbf{H}=(Df)^{\bot},
\end{equation}
where $\textbf{H}$ is the mean curvature vector of $M$ and $(Df)^{\bot}$ is the normal component of Euclidean gradient of $f$. In particular, if $f(x)=g(|x|)$, (\ref{1-3}) can be rewritten as
\begin{equation}\label{1-4}
g(|x|)\textbf{H}= \frac{g'(|x|)x^{\bot}}{|x|}.
\end{equation}
This type of variational problem is well-studied and is closely related to several problems: the mean curvature flow (cf. \cite{Hu90}, \cite{An92}, \cite{DYM10}, \cite{CM12}, \cite{AS13}, \cite{ALW14}, \cite{KM14}, \cite{LL15}, \cite{D20}), the submanifold geometry properties of weighted minimal surfaces (cf. \cite{Liu13}, \cite{M16}, \cite{CMZ14}, \cite{CMZ15}, \cite{W17}) and the isoperimetric problems (cf. \cite{RCBM08}, \cite{FMP11}, \cite{M16}, \cite{CRS16}, \cite{Cs18}, \cite{C19}, \cite{MN22}). In Riemannian $G$-manifolds, searching for closed, minimal or constant mean curvature (CMC) submanifolds which are invariant under isometry subgroups needs to find minimal or CMC submanifolds in the orbit spaces with weighted metrics (\cite{HsL71}, \cite{L72}, \cite{Hs82}, \cite{Hs83}, \cite{An92}, \cite{KM14}, \cite{Z16}), etc.

\medskip\noindent

To study a geometric variational problem, as the theory of minimal surfaces, there are some fundamental questions: the existence of solutions to the Euler-Lagrange equation, the stability of solutions, the existence of minimizers, the regularity of minimizers and the structure of the singular sets in these minimizers. For example, in studying the singularity of mean curvature flow, the functional with density function $f(x) = e^{\epsilon\frac{|x|^2}{4}}$ is often considered. For the case $\epsilon=-1$, some regularity results were obtained by T.H. Colding and W.P. Minicozzi in \cite{CM12}. For the case $\epsilon =1$, when a minimal hypercone $\mathcal{C}$ is the boundary of a domain 
$\Omega \subset \mathbb{R}^n$, Q. Ding \cite{D20} proved that $\mathcal{C}$ is stable for the area functional if and only if it is stable for $f$-weighted area functional, and he also showed that $\mathcal{C}$ is (one-sided) area-minimizing in $\Bar{\Omega}$ if and only if $\mathcal{C}$ is $f$-minimizing in $\Bar{\Omega}$, i.e. minimizing in one of the components of $\mathbb{R}^n\setminus\mathcal{C}$, Ding also proved that graphic self-expanding hypersurfaces are $f$-minimizing. The Euler-Dierkes-Huisken variational problem is related to the famous Bernstein problem and Hilbert's 19th problem which concerns the regularity of minimizers of variation integrals (see a nice survey by C. Mooney \cite{Mo22}). In this paper, we mainly study the stability of the $f$-stationary solutions and minimizers of the Euler-Dierkes-Huisken functional.

\medskip\noindent
\medskip\noindent

\subsection{Definitions of regular, non-regular and $f$-minimizing cones}

\medskip\noindent

In the context of geometric measure theory (GMT), area-minimizing surfaces in $\mathbb{R}^{n}$ are integral currents (cf. \cite{FF60}) that globally minimize the area functional with a given boundary. An important feature is that every interior singular point $p$ of an area-minimizing surface $S$ has at least one oriented tangent cone as the blowing-up limit of $S$, this cone is itself area-minimizing (cf. 4.3.16 and 5.4.3 in \cite{Fe69}). Therefore, studying area-minimizing cones helps in understanding the structures of singularities of area-minimizing surfaces.
In this subsection, we will introduce the concepts of regular and non-regular cones, the space of rectifiable currents reduced modulo $p$, and the definition for minimizing cones with respect to radial density
functions in the Euler-Dierkes-Huisken functional.

\medskip\noindent

\medskip\noindent

\textbf{Definition 1.1.} \textit{For a submanifold $($or rectifiable subset$)$ $ \Sigma^{k-1} \subset S^{n-1}(1)$, the cone over $\Sigma$ is defined by
$$
\mathcal{C}=C(\Sigma):=\{ t x\;| t \geq 0, x \in \Sigma \} \subset \mathbb{R}^{n},
$$
and $\Sigma$ is called the \textit{link} of $\mathcal{C}$. We also define $\mathcal{C}_{1}$ to be the truncated part of $\mathcal{C}$ inside the unit ball, i.e., $\mathcal{C}_{1}:=\mathcal{C} \cap \textbf{B}^{n}(1)$, and so $\partial \mathcal{C}_1= \Sigma$.}

\medskip\noindent

Note that the origin $O$ becomes a natural singularity of $\mathcal{C}$ when $\mathcal{C}$ is nonplanar. If $\Sigma$ is a closed smooth submanifold in unit spheres, then $\mathcal{C}$ only carries an isolated singularity at the origin. These cones are referred to as \textit{regular cones}, which were first defined in \cite{HS85} for cones over closed smooth hypersurfaces. And any other cones with larger singular sets will be called \textit{non-regular cones} here. 

\medskip\noindent

As pointed out by Lawlor in \cite{La91} and \cite{La98}, the methods of curvature criterion and directed slicing both can deal with the cones (surfaces) that are unorientable. In such cases, the space of surfaces is often taken to be the \textit{rectifiable currents reduced modulo 2}.
\medskip\noindent

\textbf{Definition 1.2.} (see \cite{La98},\cite{KL99}, \cite[Section 11.1]{M16}) \textit{Let $p\geq$ 2 be a integer, two rectifiable currents $T$, $S$ are congruent modulo $p$ if $T-S=p Q$ for some rectifiable current $Q$, say $T\equiv S,\mod p$. We define rectifiable currents reduced modulo $p$ to be the congruent classes modulo $p$.}

\medskip\noindent

This definition is a simplified version of the original one in \cite[4.2.26]{Fe69}, which involves the process of taking flat norm closure. For their relation, Morgan (\cite[Section 11.1]{M16}) gave a nice talk, for a detailed study of area-minimizing currents mod $p$, one can refer to \cite{DHMS20}. It is worth noting that, even though the existence theory of minimizers is not applicable for currents reduced modulo $p$, it does not pose any problems since we generally start with a proposed minimizer, as explained in \cite[4.2.1]{La98}.
 
\medskip\noindent

\textbf{Definition 1.3.}  \textit{Let $\mathcal{C}=C(\Sigma)$ be a $k$-dimensional cone in Euclidean space which can either be regular or non-regular. If the truncated cone $\mathcal{C}_1$ is orientable, we consider it as a rectifiable current, and if $\mathcal{C}_1$ is unorientable, we consider it as a rectifiable current reduced mod 2. We say that the cone $\mathcal{C}$ is area-minimizing if
$$
\mathcal{H}_k(\mathcal{C}_{1}) \leq \mathcal{H}_k(S)
$$
for any rectifiable current $($rectifiable current reduced mod 2$)$ $S$ which satisfies $\partial S=\Sigma$ $(\partial S\equiv \Sigma, \mod 2)$.
Moreover, we call $\mathcal{C}$ is $f$-minimizing if
$$
\mathcal{E}_f(\mathcal{C}_{1}) \leq \mathcal{E}_f(S)
$$
for any rectifiable current $($rectifiable current reduced mod 2$)$ $S$ which satisfy $\partial S=\Sigma$ $(\partial S\equiv \Sigma, \mod 2)$.}

\medskip\noindent

A direct calculation shows that $\mathcal{H}_{k-1}(\Sigma)<\infty$ implies $\mathcal{H}_k(\mathcal{C}_{1})< \infty$, and for a given radial density function $f(x)=g(|x|)$, $\mathcal{E}_f(\mathcal{C}_{1})< \infty$  if and only 
if
\begin{equation}\label{priori}
\int_{0}^1g(t) t^{k-1} dt < \infty,
\end{equation}
i.e. the Plateau's problem of minimizing cone for Euler-Dierkes-Huisken functional is valid if \eqref{priori} is satisfied, which is \textit{a priori condition}. Similarly, for the study of $f$-stable minimal cone, we assume \eqref{priori} is satisfied. For instance, in Euler-Dierkes-Huisken functional with 
$f(x)=|x|^{\alpha}$, it is necessary that $\alpha >-k$. It is clear that an oriented $f$-minimizing cone is naturally $f$-stable since the comparison surfaces include all small deformations of $\mathcal{C}_1$.

\medskip\noindent

It is necessary to introduce some typical examples of area-minimizing cones, including both regular and non-regular ones.
The study of area-minimizing regular cones dates back to the famous \emph{Simons cones} ($p=q$) and \textit{Lawson cones}(\cite{L72})
$$
\mathcal{C}_{p,q}:=C(S^p \times S^q)=\{(x_{0},\ldots,x_{p+q+1})\in \mathbb{R}^{p+q+2}|\;q(x_{0}^2+\cdots+x_{p}^2)=p(x_{p+1}^2+\cdots+x_{p+q+1}^2) \},
$$ which plays important roles in studying Bernstein problem (cf. \cite{Sim68}, \cite{BDGG69}). It is well-known that $\mathcal{C}_{p,q}$ is area-minimizing if and only if $p+q\geq 7$ or $p+q=6,|p-q|\leq 3$ (cf. \cite{Sim68}, \cite{BDGG69}, \cite{L72}, \cite{Si74}, \cite{La91}, \cite{Da04}, \cite{DP09}, etc). In 1991, Lawlor \cite{La91} provided a general method, called \emph{curvature criterion}, to prove a given regular cone's area-minimization. This method is powerful and has led to the discovery of many new area-minimizing regular cones, as in \cite{La91}, \cite{Ke94}, \cite{Ka02}, \cite{XYZ18}, \cite{TZ20}, \cite{JC22}, \cite{JCX22} etc. For more historical notes and methods of proof concerning a given regular minimal cone's area-minimization, one can refer to Section 2 in \cite{La98} for details. 

\medskip\noindent

It is a challenging task to search for area-minimizing non-regular cones.  Canonical examples of such cones include singular complex projective varieties when seen as cones in Euclidean spaces, and the union of two oriented $k$-planes under a characteristic angle condition (cf. \cite{M83}, \cite{Ma87}, \cite{La89}), etc. After the work \cite{La91}, Lawlor \cite{La98} developed another powerful but more general method, called \emph{directed slicing}, to study highly-singular minimizing surfaces. Roughly speaking, by using the area-coarea formula, it turns the area-minimizing problem into a family of localized weighted Plateau's problems. This also motivates us to study minimizing cones associated with weighted area integrals. For examples of non-regular area-minimizing cones, readers are referred to \cite{La98}, \cite{LM96}, \cite{KL99}, \cite{CJX24}, etc.

\medskip\noindent
\medskip\noindent

\subsection{Main results and organizations}

\medskip\noindent

In Section 2, we first derive the Euler-Lagrange equation of Euler-Dierkes-Huisken functional and then generalize Dierkes-Huisken's results (\cite{DH23}) for Euler-Dierkes-Huisken functional with density function $f(x)=|x|^{\alpha}$ to high codimensional cases with a more general radial non-negative density function $f(x)=g(|x|)$, see Theorem 2.3 for details.

\medskip\noindent

We first investigate the $f$-stability of $f$-stationary surfaces. For the planes, we have

\medskip\noindent

\textbf{Theorem 1.4.} \textit{The planes $\mathbb{R}^k$ through origin are $f$-stable for the Euler-Dierkes-Huisken functional with density function $f(x)=g(|x|)$, where $g(t)\geq 0$ and $g'(t)\geq 0$}.

\medskip\noindent

However, as a corollary of \cite[Theorem 9.2]{CM12}, the hyperplanes through origins are not $f$-stable with $f(x)={\rm e}^{-\frac{|x|^2}{4}}$. The following theorem says that the high codimensional generalization of (ii)
in \cite[Theorem 1.6]{DH23} is not true.

\medskip\noindent

\textbf{Theorem 1.5.} \textit{The $k$-dimensional minimal submanifolds in $S^{n-1}(R)(R>0)$ as submanifolds in $\mathbb{R}^{n}$ are $f$-unstable, where $f(x)=|x|^{-k}$ and $1\leq k \leq n-2$.}

\medskip\noindent

Next, we consider the density function $f(x)$ takes the form $f(x)={\rm e}^{h(|x|)}$, $h(t)$ is $C^2$ on $\mathbb{R}^+$. To find $f$-stable surfaces, 
we need to consider some further assumptions on $h(t)$, i.e.,

(A-1) \hspace{0.2cm} $\lim \limits_{t \rightarrow 0^+} {\rm e}^{h(t)} t^{k}=a<+\infty$;

(A-2) \hspace{0.2cm} $\inf\limits_{t\in(0,+\infty)}th'(t)=b>2-k$;

(A-3) \hspace{0.2cm} $t(h')^2+2th''+2nh'\geq 0$ for $t>0$;

(A-4) \hspace{0.2cm} $\int_{0}^r {\rm e}^{h(t)} t^{k-1} dt < +\infty$, for any $r>0$.

\noindent When $h(t)$ is non-decreasing, (A-1) is naturally satisfied and (A-2) is satisfied for $k\geq 3$, also it is worth pointing out that the assumption (A-1) and (A-2) implies that (A-4) holds. For the minimal submanifolds
in unit spheres with globally flat normal bundles (cf. \cite{T87} and \cite{PT88}), we have

\medskip\noindent

\textbf{Theorem 1.6.} \textit{Let $\Sigma$ be a closed minimal $(k-1)$-dimensional submanifold in the unit sphere $S^{n-1}(1)$ with the globally flat normal bundle, if the squared norm of the second fundamental forms of $ \Sigma$ satisfies $|\textbf{A}_{\Sigma}|^2 \leq \frac{(k-2+b)^2}{4}+b$, then the cone over $\Sigma$ is $f$-stable, where $f(x)={\rm e}^{h(|x|)}$, $h(t)$ satisfies} (A-1), (A-2).

\medskip\noindent

If (A-1) is satisfied and $b=0$ in (A-2), this theorem is closely related to a well-known fact that: for an oriented regular minimal hypercone $C(\Sigma)$ in $\mathbb{R}^n$ if pointwisely $|\textbf{A}_{\Sigma}|^2 \leq \frac{(n-3)^2}{4}$, then $C(\Sigma)$ is a stable minimal cone for the area functional (cf. \cite[Appendix B]{Simon83}). For the Euler-Dierkes-Huisken functional, $h(t)=\alpha \ {\rm ln} t$ and $b=\alpha$ if the density function $f(x)=|x|^\alpha$, $h(t)=\frac{t^2}{4}$ and $a=0$ for self-expander. However, the assumption (A-2) is not satisfied for self-shrinkers. 

\medskip\noindent 

The submanifolds in unit spheres having globally flat normal bundles will span minimal cones with globally flat
normal bundles. These submanifolds include the oriented minimal hypersurfaces, the principle orbits of isometric polar group actions (for example $s$-representation, see \cite[Section 2.3.2]{BCO16} and \cite[Theorem 5.7.1]{PT88}), and the product of spheres. The product of two or more spheres has a property that: $|\textbf{A}_{\Sigma}|^2={\rm dim}\ \Sigma=k-1$ (cf. \cite[Theorem 5.1.1]{La91}). Therefore, from Theorem 1.6, we have a corollary that generalizes the result of Dierkes-Huisken for cones over the product of two spheres (cf. \cite[Corollary 2.2]{DH23}).

\medskip\noindent

\textbf{Corollary 1.7}. \textit{The $k$-dimensional minimal cones $\mathcal{C}=C \left(S^{k_1} \times \cdots \times S^{k_{l}}\right)$ are $f$-stable, where $f(x)=|x|^{\alpha}$
with $\alpha \geq -k+2\sqrt{2(k-1)}$. Particularly, these minimal cones are $f$-stable if ${\rm dim} \ \mathcal{C}=k \geq 7$ and $\alpha \geq 4\sqrt{3}-7 $.}

\medskip\noindent

There is a family of minimal submanifolds in the unit spheres with globally flat normal bundles constructed in \cite[Lemma 4.4.4]{La91}, which can be seen as a generalization for the product of spheres. A similar construction, called \emph{minimal product}, which has been studied in \cite{X03}, \cite{CH18}, \cite{TZ20}. Their geometric properties and area-minimizing problems are also studied in \cite{TZ20}, \cite{JCX22}, \cite{WW23}, etc. More explicitly, let $\Sigma_i \subset S^{n_i-1}(r_i)$, $1\leq i \leq m$, be minimal submanifolds with $\dim \Sigma_i=k_i$, whose normal bundles are globally flat and the second fundamental forms satisfy $\sup\limits_{|v_i|=1} |\textbf{A}_{\Sigma_i}^{v_i}|^2 \leq k-1$. Here, $k-1=\sum\limits_{i=1}^{m}k_i$ and $r_i^2= \frac{k_i}{k-1}$. So,
$\Sigma=\Sigma_1 \times \cdots \times \Sigma_m$ is minimal in the unit sphere $S^{n-1}(1)$, where $n=\sum\limits_{i=1}^{m}n_i+m$. Then, the cone $C \left(\Sigma_1 \times \cdots \times \Sigma_m\right)$ are 
$f$-stable, where the density function $f(x)=|x|^{\alpha}$ and $\alpha \geq -k+2\sqrt{2(k-1)}$. This is also a generalization of Corollary 1.7.

\medskip\noindent

The next result is inspired by the work of Q. Ding \cite{D20}.

\medskip\noindent

\textbf{Theorem 1.8}. \textit{The oriented stable minimal hypercone in $\mathbb{R}^n$ is also $f$-stable, where $f(x)={\rm e}^{h(|x|)}$ and $h(t)$ satisfies} (A-3), (A-4).

\medskip\noindent

The term $t(h')^2+2th''+2nh'$ in (A-3) comes from the difference between the $f$-stability inequality for $\xi$ and the stability inequality of area-functional for $\xi {\rm e}^{\frac{h}{2}}$, where $\xi$ is some compactly supported function on the cone. For the density function $f(x)=|x|^\alpha$, we have $h'(t)= \frac{\alpha}{t}$
and $h''(t)= -\frac{\alpha}{t^2}$, then the assumption (A-3) and (A-4) are satisfied if $\alpha\geq 0$.
So, We have

\medskip\noindent

\textbf{Corollary 1.9}. \textit{The oriented stable minimal hypercone in $\mathbb{R}^n$ is also $f$-stable, where $f(x)= |x|^\alpha$ and $\alpha \geq 0$.}

\medskip\noindent
\medskip\noindent

In Section 3, we recall the methods of Lawlor curvature criterion and Lawlor directed slicing, which will be used to prove the $f$-minimizing property of higher codimensional, regular or non-regular cones in Section 4. These methods are different from Bombieri-De Giorgi-Giusti's method and the theory of sets of finite perimeters used in \cite{DH23}, which are 
only suitable for the minimal hypercones.

\medskip\noindent

For the regular minimal cones, we have

\medskip\noindent

\textbf{Theorem 1.10}.  \textit{The regular minimal cones satisfy Lawlor curvature criterion are $f$-minimizing
in the sense of Definition 1.3, where $f(x)=g(|x|)$ is non-negative and $g$ is non-decreasing. Also, the planes through the origin are $f$-minimizing for these integrals.} 

\medskip\noindent

The proofs of Theorem 1.10 and Theorem 1.12 below are based on the area-coarea formula in GMT (cf.\cite[3.2.22]{Fe69} and \cite[Theorem 4.5.1]{La98}). So, the density function only needs to be integral in the Hausdorff measure
when considering the $f$-minimizing property of surfaces. The radially non-decreasing density functions in Theorem 1.10 and Theorem 1.12 naturally satisfy the a priori condition \eqref{priori}, which is necessary in considering the $f$-minimizing property. 
 
\medskip\noindent 

Minimal cones satisfying Lawlor curvature criterion are actually area-minimizing, and there are plenty of minimal cones that satisfy the Lawlor curvature criterion, see Remark 3.6. It is worth mentioning that the statement of planes
is a natural generalization of Dierkes-Huisken's work (cf. \cite[Theorem 3.3]{DH23}). On the other hand, inspiring
from \cite[Theorem 5.1.1]{La91}, we have the following corollary of Theorem 1.10. That is

\medskip\noindent 

\textbf{Corollary 1.11}. \textit{The minimal cone $\mathcal{C}=C \left(S^{k_1} \times \cdots \times S^{k_{m}}\right)$ are $f$-minimizing with $\dim \mathcal{C} \geq 7$ and $k_i>1$, where $f(x)=g(|x|)$ and $g$ is non-decreasing. Particularly, such minimal cones are $|x|^{\alpha}$-minimizing for $\alpha \geq 0$}. 

\medskip\noindent

For the non-regular cones, we recall the definitions of determinantal varieties and Pfaffian varieties 
(skew-symmetric determinantal varieties), one can refer to \cite{CJX24} for more details.
Let $M(p,q;\mathbb{R})$, $p\leq q$, be the space of matrices of size $p\times q$ over $\mathbb{R}$, and let $Skew(m,\mathbb{R})$ be the space of $m \times m$ skew-symmetric matrices over $\mathbb{R}$. The \textit{determinantal variety} $C(p,q,r)$ with $r< p$ are defined by 
$$
C(p,q,r):=\{X \in M(p,q;\mathbb{R}): {\rm rank} \ X \leq r \}.
$$
The \textit{Pfaffian variety} ${\bf C}(m,2r)$ with $2r<m$ are defined by
$$
{\bf C}(m,2r)=\{X \in Skew(m,\mathbb{R}): {\rm rank} \ X \leq 2r \}.
$$
Those minimal varieties have stratified singular sets. By using Lawlor directed slicing, the area-minimizing properties of cones $C(p,q,r)$, ${\bf C}(m,2r)$ have been studied in \cite{KL99} and \cite{CJX24}, respectively. In this paper, we will study the $f$-minimizing properties of these two kinds of varieties with respect to the Euler-Dierkes-Huisken functional. That is

\medskip\noindent

\textbf{Theorem 1.12}.  \textit{\emph{(1)} The oriented area-minimizing cones $C(p,q,r)$ with $p+q-2r\geq 4$ and $p+q$ is even, ${\bf C}(m,2r)$ with $m-2r\geq 3$, are $f$-minimizing, where $f(x)=g(|x|)$ is non-negative and $g$ is non-decreasing.}

\medskip\noindent 

\textit{\emph{(2)} The unoriented minimal cones $\mathcal{C}:=C(p,q,r)$ with $p+q-2r\geq 4$ and $p+q$ is odd, are $f$-minimizing in the sense of Definition 1.3, where $f(x)=g(|x|)$ is non-negative and $g$ is non-decreasing.}

\medskip\noindent

\textit{\emph{(3)} The minimal cones $C(2,2,1)=C(S^1 \times S^1)$, $C(2,3,1)$ and ${\bf C}(4,2)=C(S^2 \times S^2)$ are $f$-minimizing in the sense of Definition 1.3, where $f(x)=g(|x|)$ is non-negative and $\frac{g(t)}{t}$ is non-decreasing, such as: $|x|^{\alpha}(\alpha \geq 1)$.}

\medskip\noindent

This theorem generalizes related results in \cite{KL99},\cite{CJX24}. Particularly, $C(2,3,1)$ belongs to a family of unsolved cases for the area-minimization: $C(p,p+1,p-1)$ in \cite{KL99}, and it is shown $f$-minimizing here.  Also, notice that the density function $|x|^\alpha$, $\alpha\geq 1$, satisfies the condition in (3) of Theorem 1.12, so combining it with Corollary 1.11 implies the following corollary, which is a generalization of Dierkes-Huisken's work (cf. \cite[Theorem 3.1]{DH23}). 

\medskip\noindent

\textbf{Corollary 1.13}. \textit{All Simons cones $C(S^{p} \times S^{p})$, $p\geq 1$, are $|x|^{\alpha}$-minimizing for $\alpha \geq 1$.}

\medskip\noindent

Since the method of compensating for the composite weighting functions along parabolic secondary slicing sets is \textbf{false} when dealing with the generic unsolved cases\footnote{Although the set-theoretical tangent cones of determinantal hypersurfaces $C(p,p,p-1)$ and Pfaffian hypersurfaces ${\bf C}(2n,2n-2)$ are proved not to be area-minimizing at some ''nearly regular'' points in \cite{KL99} and \cite{CJX24}, we still prefer not to conclude that these varieties are definitely non-minimizing, see a discussion in \cite[Section 1.3]{CJX24}.} in \cite{KL99} and \cite{CJX24}: $C(p,p+1,p-1)$, $C(p,p,p-1)$ and ${\bf C}(2n,2n-2)$, so we don't know that
whether they are $f$-minimizing except for the three cones in Theorem 1.12 (3), see Proposition 4.1 for details. Finally, we should point out there is very little known about $f$-minimizing cones of Euler-Dierkes-Huisken variational problems for decreasing radial density functions $f(x)=g(|x|)$, like $f(x)=|x|^\alpha$, $-k<\alpha <0$, etc.

\medskip\noindent
\medskip\noindent

\section{$f$-stationary and $f$-stable submanifolds}
\subsection{$f$-stationary submanifolds}

\medskip\noindent

Let $M$ be an $k$-dimensional $C^2$ submanifold embedded in $\mathbb{R}^{n}$. $\Bar{M}$ denote the closure of $M$. Let $U \subset \mathbb{R}^{n}$ be a open subset, with $U \cap M \neq \varnothing$, $U \cap (\Bar{M}\setminus M)=\varnothing$. Suppose there is a 1-parameter family of $C^2$ maps $\varphi_t: M \rightarrow \mathbb{R}^{n}$, $-\varepsilon <t < \varepsilon$, such that: 
$$
\left\{\begin{array}{l}
\varphi_{0}(x):=\varphi(0, x)=x, \quad \forall x \in M, \\
\varphi_{t}(x):=\varphi(t, x)=x, \quad \forall(t, x) \in(-\varepsilon, \varepsilon) \times(M \backslash K),
\end{array}\right.
$$
where $K \subset U$ is a compact subset. Let 
$ X(x)=\left.\frac{\partial \varphi}{\partial t}(t, x)\right|_{t=0} $ and 
$ Z(x)=\left.\frac{\partial^2 \varphi}{\partial t^2}(t, x)\right|_{t=0}$
denote the initial velocity and acceleration vectors of $\varphi_t$ respectively. Then $X$ and $Z$ are compactly supported in $K \subset U$ and clearly 
$$
\varphi_t(x)=x+t X(x)+\frac{t^2}{2} Z(x)+o\left(t^2\right) .
$$

\medskip\noindent 

We call $\varphi_{t}$ a compactly supported deformation of $M$ in $U$, also denote $M_t=\varphi_t(M)$, for the variation integral 
$$
\mathcal{E}_f(M):=\int_M f(x) d \mathcal{H}_k,
$$
here in Section 2, we require $f: \mathbb{R}^{n} \rightarrow \mathbb{R}$ to be some function of class $C^2$, and the requiremens for regularity of $f$ and $M$ will be released in Section 4.

\medskip\noindent 

We let 
$$
\delta \mathcal{E}_f(M, \varphi_{t})=\left.\frac{d}{d t}\right|_{t=0} \int_{M_t} f(x) d \mathcal{H}_k \ \
{\rm and} \ \
\delta^2 \mathcal{E}_f(M, \varphi_{t})=\left.\frac{d^2}{d t^2}\right|_{t=0} \int_{M_t} f(x) d \mathcal{H}_k 
$$
denote the first and second variations.

\medskip\noindent 

Let  $\langle, \rangle$ denotes the Euclidean inner product (Euclidean metric), if $X$ is a vector field on $M\subset \mathbb{R}^{n}$, we let $X^{\top}$ and $X^{\bot}$ denote the tangential and normal components of $X$.
For convenience, we also choose $\{e_1,\ldots,e_k\}$ and $\{e_{k+1},\ldots,e_{n} \}$ to be the local tangential and normal orthonormal frames of $M$, respectively. The covariant derivative $D$ on $\mathbb{R}^{n}$ induces a covariant derivative (in fact, Riemannian connection) $\nabla$ on $M$ and second fundamental form $\textbf{B}$ of $M$:

$$
\nabla=(D)^{\top} \ {\rm and} \ \textbf{B}(X,Y)=(D_{X}Y)^{\bot},
$$
also denote the Hessian of $f(x)$ in $\mathbb{R}^{n}$ by ${\rm Hess} f$.

\medskip\noindent 

For $M$ at least $C^2$, we can define its second fundamental forms in $\mathbb{R}^{n}$ and the mean curvature vector (for convenience, we take the mean curvature vector to be the sum rather than the average of the trace of second fundamental forms):
$$
\textbf{H}:= \sum_{i=1}^{k} \textbf{B}(e_i,e_i) =\sum_{i=1}^{k} \left(  D_{e_{i}} e_{i} \right)^{\bot}.
$$

\medskip\noindent 

Let $M^{k} \subset \mathbb{R}^{n}$ and $f: \mathbb{R}^{n} \rightarrow \mathbb{R}$ be as above. By chain rules: 
$$
\frac{d}{d t}(f(\varphi_{t}(x)))= \langle Df (\varphi_{t}(x)), \frac{d}{d t} \varphi_{t}(x)\rangle,
$$
$$
\frac{d}{d t} \langle Df (\varphi_{t}(x)), \frac{d}{d t} \varphi_{t}(x)\rangle = \langle  {\rm Hess} f \cdot \frac{d}{d t} \varphi_{t}, \frac{d}{d t} \varphi_{t}(x)\rangle + \langle Df (\varphi_{t}(x)), \frac{d^2}{d t^2} \varphi_{t}(x)\rangle 
$$ and the standard computations for Jacobian $J_{\varphi_{t}}$ (for example, see \cite[Chapter 9]{Simon83}), then the first and second variation of (\ref{fequation}) can be given by the formulas (see Proposition 1.1 in \cite{DH23}):

\begin{equation}\label{1var}
  \delta \mathcal{E}_f(M, \varphi_{t})=\int_M\left\{f(x) \operatorname{div} X+ \langle Df, X\rangle  \right\} d \mathcal{H}_k  
\end{equation}
and 
\begin{equation}\label{2var}
\begin{aligned}
\delta^2 \mathcal{E}_f(M, \varphi_{t})= & \int_M f\left[\operatorname{div} Z+(\operatorname{div} X)^2+\sum_{i=1}^k\left|\left(D_{e_i} X\right)^{\perp}\right|^2
-\sum_{i, j=1}^k \langle e_i, D_{e_j} X \rangle \langle e_j, D_{e_i} X \rangle \right] d \mathcal{H}_n \\
&+\int_M \left[ 2 \langle Df, X\rangle \operatorname{div} X+ {\rm Hess} f(X,X)+\langle Df, Z\rangle\right] d \mathcal{H}_k
\end{aligned}
\end{equation}
where $\operatorname{div} X=\sum\limits_{i=1}^{k} \langle D_{e_i}X,e_{i} \rangle$ represents the divergence of $X$ on $M$.

\medskip\noindent

\textbf{Definition 2.1. (see Definition 1.2 in \cite{DH23})}.  \textit{A $C^1$ submanifold $M$ is called stationary (or $f$-stationary) in $U\subset \mathbb{R}^{n}$ if 
$$
\delta \mathcal{E}_f(M, \varphi_{t})=0
$$
for any compactly supported deformation $\varphi_{t}$ in $U$.}

\medskip\noindent 

The first variation formula was well-studied:

\medskip\noindent 

\textbf{Proposition 2.2. (see \cite{I95}, \cite{Ec04}, \cite{RCBM08}, \cite{DH23})}. \textit{Suppose $M$ is of class $C^2$, then $M$ is $f$-stationary in $U$ if and only if the equation \eqref{1-3} is satisfied. For $f(x)=g(|x|)$, equation \eqref{1-3} is equivalent to equation \eqref{1-4}.}

\medskip\noindent 

Searching for the general solutions of \eqref{1-4} is not easy, for self-shrinking solutions in MCF, there are plenty of complex examples such as \textit{Angenent torus} \cite{An92}, etc. However, easy to see regular minimal cones or planes through the origin are solutions.

\medskip\noindent 

For variation integral: $f(x)=|x|^{\alpha}$, we have the following  theorem:

\medskip\noindent

\textbf{Theorem 2.3 (compare with \cite[Theorem 1.6]{DH23})}. \textit{For $\alpha \in \mathbb{R}$ and $f(x)=|x|^\alpha$, let $M$ be a $k$-dimensional, $f$-stationary submanifold in $\mathbb{R}^{n}$ or $\mathbb{R}^{n}\setminus\{0\}$, i.e. }
\begin{equation}\label{mod}
{\rm \textbf{H}}=\alpha |x|^{-2} x^{\bot} 
\end{equation}
\textit{then we have:}

\medskip\noindent 

\textit{\emph{(1)} If $\alpha>-k, M$ cannot be compact;}

\medskip\noindent 

\textit{\emph{(2)} If $\alpha=-k$, then as submanifolds in $\mathbb{R}^{n}$, all $k$-dimensional minimal submanifolds in the hyperspheres $S^{n-1}(R)$, $R>0$ satisfy \eqref{mod};}

\medskip\noindent 

\textit{\emph{(3)} If $\alpha<-k$ and $M \subset \mathbb{R}^{n}\setminus\{0\}$ of class $C^2$ is $f$-stationary in $\mathbb{R}^{n}\setminus\{0\}$, then its closure $\Bar{M}$ contains the origin.}

\medskip\noindent 
\medskip\noindent 

\textbf{Proof}:
The original proof in \cite{DH23} can be modified and applied to general cases.

\medskip\noindent 

(1) By a standard computation (see \cite{CM11}), for a fixed vector $a$ (also denoted its parallel vector field in $\mathbb{R}^{n}$)
$$
e_i \langle x, a \rangle= \langle D_{e_i}x, a \rangle = \langle e_i, a \rangle 
$$
since $x$ is the position vector field, then
$$
\Delta \langle x, a \rangle=e_i e_i \langle x, a \rangle- \nabla_{e_i}e_{i} \langle x, a \rangle= e_{i} \langle e_{i}, a \rangle-\langle \nabla_{e_i}e_{i}, a \rangle=\langle D_{e_i}e_{i}-\nabla_{e_i}e_{i}, a \rangle=\langle \textbf{H}, a \rangle
$$
yields $\Delta x= \textbf{H}$, where we have used Einstein Summation for simplifying expressions.

Choose a standard orthonormal basis of $\mathbb{R}^n$ by $\{E_ A\}(1\leq A\leq n)$, then,
$$
\begin{aligned}
\Delta |x|^2 & =\sum_{A,i} e_i e_i \langle x, E_{A} \rangle^2  - \sum_{A,i} \nabla_{e_i}e_{i} \langle x, E_{A} \rangle^2 \\
&= 2 \sum_{A} \langle x, E_{A} \rangle \Delta \langle x, E_{A} \rangle + 2 \sum_{A,i} \langle e_i, E_{A}\rangle ^2 \\
&= 2 \sum_{i=1}^k |e_{i} |^2 + 2 \langle x,  \Delta x \rangle \\
& =2 k+2 \alpha \frac{|x^{\bot}|^2}{|x|^2} > \frac{2k|x^\top|^2}{|x|^2} \geq 0
\end{aligned}
$$

Since $M$ is compact, by Hopf's maximum principle we must have $|x|^2=r^2$ for some constant $r>0$, thus $M$ is contained in a sphere of radius $r$ which centered at the origin, so $x$ is naturally a normal vector field. Then \eqref{mod} is equivalent to $\textbf{H}=\frac{\alpha}{r^{2}} x$, it follows that $\langle \textbf{H}, x\rangle = \alpha$, and this leads to a contradiction since a direct computation shows that $\langle \textbf{H}, x\rangle = -k$. 

\medskip\noindent 

(2) A direct corollary of Takahashi's theorem: $M$ is a $k$-dimensional minimal submanifolds in $S^{n-1}(R)$ iff $\Delta x + \frac{k}{R^2}x=0$.

\medskip\noindent 

(3) Conversely, assume that $ d:= inf_{x\in M} |x| > 0$, then there exists a sequence $x_j \in M$ with $|x_j| \rightarrow d$, and also $ x_j \rightarrow x_0 \in \Bar{M}$, as $j \rightarrow \infty $. 
By assumption, we have $x_0 \in M$ and  $|x_0|=d=inf_{x \in M}|x|>0$. It follows that $\nabla |x|^2(x_0)=2x^{\top}(x_{0})=0$, i.e. $x_0^{\top}=0$, then 
$$
\Delta |x|^2 =2 k+2 \alpha \frac{|x^{\bot}|^2}{|x|^2} < \frac{2k|x^\top|^2}{|x|^2}, 
$$
the last expression equals zero at $x_0$ which is away from the origin at a distance $d>0$, thus there exists a neighborhood $U_{\epsilon}(x_0) \cap M $, such that $\Delta |x|^2<0$ on it, by Hopf's maximum principle, it follows that $|x|^2=r^2$ for some constant $r>0$ on $U_{\epsilon}(x_0) \cap M$, by a similar discussion in the proof of (1), this will also lead to a contradiction. $\Box$

\medskip\noindent 

\textbf{Remark 2.4}. \textit{Recall the formula \eqref{priori},  when we study the $f$-stable and $f$-minimizing properties of a minimal cone $\mathcal{C}$ for variation integral: $f(x)=|x|^{\alpha}$, to ensure the local boundedness of $\mathcal{E}_f(\mathcal{C})$, we always assume $\alpha>-k$, $k={\rm dim} \ \mathcal{C}$.}

\medskip\noindent 
\medskip\noindent 

\subsection{$f$-stable submanifolds}

\medskip\noindent

In this subsection, we assume that the a priori condition \eqref{priori} is always satisfied, and investigate the stability of $f$-stationary submanifolds in $\mathbb{R}^{n}$ or $\mathbb{R}^{n}\setminus\{0\}$. If $ M$ is $f$-stationary in $U \subset \mathbb{R}^{n} $, then \eqref{2var} reduces to (independent with $Z$)

\begin{equation}\label{rvar}
\begin{aligned}
\delta^2 \mathcal{E}_f(M, \varphi_{t})= & \int_M f\left[(\operatorname{div} X)^2+\sum_{i=1}^k\left|\left(D_{e_i} X\right)^{\perp}\right|^2-\sum_{i, j=1}^k \langle e_i, D_{e_j} X \rangle \langle e_j, D_{e_i} X \rangle \right] d \mathcal{H}_k \\
& + \int_M\left\{2 \langle Df, X\rangle \operatorname{div} X +{\rm Hess} f(X,X) \right\} d \mathcal{H}_k.
\end{aligned}
\end{equation}

\medskip\noindent 

We also assume the variation field $X$ is normal to $M$, then $\operatorname{div} X= - \langle \textbf{H}, X \rangle $,  $\langle Df, X \rangle=f \langle \textbf{H}, X \rangle$ and $\langle e_i, D_{e_j} X \rangle \langle e_j, D_{e_i} X \rangle=\langle \textbf{B}(e_i,e_j),X \rangle ^2 $. Thus \eqref{rvar} reduces to

\begin{equation}\notag
\delta^2 \mathcal{E}_f(M, \varphi_{t})= \int_M f\left[\sum_{i=1}^k\left|\left(D_{e_i} X\right)^{\perp}\right|^2- \sum_{i, j=1}^k \langle \textbf{B}(e_i,e_j),X \rangle ^2 - \langle \textbf{H}, X \rangle ^2 \right] + {\rm Hess} f(X,X) d \mathcal{H}_k.
\end{equation}

\textbf{Definition 2.5. (cf. Definition 6 in \cite{CMZ15} and Definition 1.5 in \cite{DH23})} \textit{A $C^2$ submanifold $M$ which is $f$-stationary in $U\subset \mathbb{R}^{n}$ is called $f$-stable in $U$ iff 
\begin{equation}\label{r2var}
 \int_M f\left[\sum_{i, j=1}^k \langle \textbf{B}(e_i,e_j),X \rangle ^2 + \langle \textbf{H}, X \rangle ^2 \right] d \mathcal{H}_k \leq  \int_M f \sum_{i=1}^k \left|\left(D_{e_i} X\right)^{\perp}\right|^2 + {\rm Hess} f(X,X) d \mathcal{H}_k
\end{equation}
for any compactly supported normal variation field $X$ in $U$.}

\medskip\noindent

For positive integrals $f(x)={\rm e}^{h(x)}$, \eqref{r2var} writes as (also see the Appendix in \cite{CMZ15}):

\begin{equation}\label{cmz15}
\int_M {\rm e}^{h} \left( \sum_{i=1}^k \left|\left(D_{e_i} X\right)^{\perp}\right|^2 + {\rm Hess} h(X,X) - \sum_{i, j=1}^k \langle \textbf{B}(e_i,e_j),X \rangle ^2 \right) d \mathcal{H}_k \geq 0
\end{equation}

\medskip\noindent

Plenty of important results for the $f$-stable hypersurfaces were obtained in \cite{CM12}, \cite{CMZ15}, \cite{D20}, \cite{DH23}, etc. 

\medskip\noindent

For radial integrals: $f(x)=g(|x|)=g(r)$
\begin{equation}\label{r5var}
{\rm Hess} f(X,X)=\frac{g'(r)}{r}|X|^2+ \frac{rg''(r)-g'(r)}{r^3} \langle x, X \rangle ^2,
\end{equation}
it follows that for $f$-stationary submanifolds which are also minimal in $\mathbb{R}^n$ or $\mathbb{R}^n\setminus\{0\}$ ($\textbf{H}=0$ implies that $x$ is a tangent vector field), \eqref{r2var} is equivalent to:

\begin{equation}\label{r4var}
 \int_M g(r) \sum_{i, j=1}^k \langle \textbf{B}(e_i,e_j),X \rangle ^2   d \mathcal{H}_k \leq  \int_M g(r) \sum_{i=1}^k \left|\left(D_{e_i} X\right)^{\perp}\right|^2 + \frac{g'(r)}{r}|X|^2 d \mathcal{H}_k
\end{equation}

Since the planes through the origin are totally geodesic: $\textbf{B}(e_i,e_j)=0$ for all $i,j$, we have proved that 

\medskip\noindent

\textbf {Theorem 2.6.(Theorem 1.4 in Introduction).} \textit{The planes $\mathbb{R}^k$ through origin are $f$-stable for the Euler-Dierkes-Huisken functional with density function $f(x)=g(|x|)$, where $g(t)\geq 0$ and $g'(t)\geq 0$ on $\mathbb{R}^+$}.

\medskip\noindent

In Section 4, we can prove that planes are $f$-minimizing under the above (weaker on the regularity) conditions.

\medskip\noindent

We now investigate the stability of $k$-dimensional ($k\leq n-2$) minimal submanifolds in hyperspheres $S^{n-1}(R)$ (as submanifolds in $\mathbb{R}^n$) for the integrals $|x|^{-k}$ (see Theorem 2.3 (2) here and \cite[Theorem 1.6]{DH23}).

\medskip\noindent

Substituting \eqref{r5var} for $g(t)=t^{\alpha}$ into \eqref{r2var} and note that $\textbf{H}=\frac{-k}{r^2}x^{\bot}=\frac{-k}{R^2}x$, \eqref{r2var} is reduced to:
\begin{equation}\label{r6var}
 \int_M  \sum_{i, j=1}^k \langle \textbf{B}(e_i,e_j),X \rangle ^2 -\frac{2}{k} \langle \textbf{H}, X \rangle ^2 + \frac{1}{k} |\textbf{H}|^2 |X|^2  d \mathcal{H}_k \leq  \int_M  \sum_{i=1}^k \left|\left(D_{e_i} X\right)^{\perp}\right|^2  d \mathcal{H}_k
\end{equation}

\medskip\noindent

Suppose now $X=\xi x$, $\xi \in C_{c}^1(U, \mathbb{R})$, then the inside integral of the left side of \eqref{r6var} reduces to 
$$
\left(\sum_{i, j=1}^k \langle \textbf{B}(e_i,e_j),x \rangle ^2-k \right) \xi^2=0,
$$
since $\langle \textbf{B}(e_i,e_j),x \rangle=\langle D_{e_i}e_{j},x \rangle=- \delta_{ij}$, meanwhile the right side of \eqref{r6var} equals $\int_M  R^2  |\nabla \xi|^2  d \mathcal{H}_k$, \eqref{r6var} is satisfied.

\medskip\noindent

Since $k\geq n-2$, for all the spherical normal variation fields: $X$, $\langle \textbf{H}, X \rangle=0$, then \eqref{r6var} is equivalent to 

\begin{equation}\label{r7var}
 \int_M  \sum_{i, j=1}^k \langle \textbf{B}(e_i,e_j),X \rangle ^2 + \frac{k}{R^2} |X|^2  d \mathcal{H}_n \leq  \int_M  \sum_{i=1}^k \left|\left(D_{e_i} X\right)^{\perp}\right|^2  d \mathcal{H}_n
\end{equation}

This is just the classical second variation for the area functional of minimal submanifolds in the ambient space: spheres of radius $R$. It is well-known that this stability inequality is always negative (cf. \cite{Sim68}), and it seems that the index of $f$-stability operator is also computed in the context of \cite[Appendix]{CMZ15}.

\medskip\noindent

So, we have proved:

\medskip\noindent

\textbf{Theorem 2.7 (Theorem 1.5 in Introduction).} \textit{For the integrals $f(x)=|x|^{-k}$, $1\leq k \leq n-2$, all the $k$-dimensional minimal submanifolds in $S^{n-1}(R) \subset \mathbb{R}^{n}(R>0)$ are unstable}.

\medskip\noindent

For general minimal cones $\mathcal{C}=C(\Sigma)$ of high codimensions, \eqref{r4var} is not easy to detect. Compared with the stability equality \eqref{cmz16} for oriented hypersurfaces, if we can decompose the terms  $\sum\limits_{i=1}^k  \left|\left(D_{e_i} X\right)^{\perp}\right|^2$ in \eqref{r4var} into a sum with disjoint normal indices, we may reduce the inequality \eqref{r4var}. And what we need are the 
submanifolds in spheres with \textit{globally flat normal bundle},  the word 'globally' ensures our integral is valid: 

\medskip\noindent

\textbf{Definition 2.8 (cf. \cite[Section 2.1]{PT88})}. \textit{Assume $M$ is an isometrically immersed submanifold in a Riemannian manifold $(N,g)$ with constant sectional curvature. Denote the normal bundle of $M$ by $T^{\perp}M$, the induced normal connection is denoted by $\nabla^{\perp}$. We call $T^{\perp}M$ flat if, for any point $x \in M$, there exists a neighborhood $U$ of $x$ and a parallel normal frame field on $U$. We call $T^{\perp}M$ globally flat if there exists a global parallel normal frame on $M$.}

\medskip\noindent

A global parallel normal frame $e_{\alpha}(k \leq \alpha \leq n-1) $ on $\Sigma \subset S^{n-1}(1)$ can be globally extended to $\mathcal{C}=C(\Sigma) \subset \mathbb{R}^n$ (away form origin) by directly setting $e_{i}(tx):=e_{i}(x), e_{\alpha}(tx):=e_{\alpha}(x), x\in \Sigma$. There are some examples of submanifolds in spheres or Euclidean spaces with globally flat normal bundles: oriented hypersurfaces; principle orbits of isometric polar group actions (cf. \cite[Theorem 5.7.1]{PT88}, \cite[Lemma 4.4.3]{La91}); isoparametric submanifolds in Euclidean spaces (cf. \cite[Theorem 6.4.5]{PT88}), etc.  

\medskip\noindent

For a minimal submanifold in the sphere, a global parallel normal vector field is called \textit{conservative} normal vector field in \cite[Section 4]{La91}. It was used to construct comparison surfaces for some minimal cone $\mathcal{C}=C(\Sigma)$ in case $\Sigma$ admits a conservative normal vector field. Based on these analyses, Lawlor finds that for some special classes of minimal cones, his curvature criterion is sufficient and necessary.  

\medskip\noindent

We now give another family of minimal submanifolds in the unit sphere with globally flat normal bundles, inspired by the construction in \cite[Lemma 4.4.4]{La91}. Independently,  it was delved into the so-called \textit{minimal product}, see \cite{X03}, \cite{CH18}, \cite{TZ20}, \cite{JCX22}, \cite{WW23}, etc.

\medskip\noindent

\textbf{Lemma 2.9}. \textit{Let $\Sigma_i \subset S^{n_i-1}(r_i)$, $1\leq i \leq m$, $\dim \Sigma_i=k_i$, with $\sum r_i^2=1$, be a family of submanifolds with globally flat normal bundles. Let
$$
\Sigma=\Sigma_1 \times \cdots \times \Sigma_m \subset S^{N-1}(1),
$$
where $N=\sum\limits_{i=1}^m n_i+m$, then $\Sigma$ has globally flat normal bundle. In special, it holds when all $\Sigma_i$ equals $S^{n_i-1}$.}

\textbf{Proof:} The arguments are the same as in the proof of \cite[Lemma 4.4.4]{La91}, we note here the codimension of $\Sigma$ is $\sum\limits_{i=1}^m (n_i-k_i)+m-1$, the term $(m-1)$ comes from a simple linear restricted equation. $\Box$

\medskip\noindent

Then, follow the standard arguments (see \cite[Appendix B]{Simon83} or the proof of \cite[Theorem 2.1]{DH23}), we can prove that 

\medskip\noindent

\textbf{Theorem 2.10 (Theorem 1.6 in Introduction)}. \textit{Assume $\Sigma$ is a closed, $C^2$ minimal $(k-1)$-dimensional submanifold in the unit sphere $S^{n-1}(1)$ with the globally flat normal bundle, if for any $t>0$, 
$th'(t)$ has a low bound $b >2-k$ and the second fundamental forms of $ \Sigma$ has a uniform upper bound $|\textbf{A}_{\Sigma}|^2 \leq \frac{(k-2+b)^2}{4}+b$, then the cone over $\Sigma$ is $f$-stable, where $f(x)={\rm e}^{h(|x|)}$.}

\medskip\noindent

\textbf{Proof:} Let $\mathcal{C}_2:=\mathcal{C}-\{0\}$ be the regular part, according to \eqref{cmz15} and \eqref{r5var}, $\mathcal{C}$ is $f$-stable iff 
\begin{equation}\label{fn1}
\delta^2 \mathcal{E}_f= \int_{\mathcal{C}_2} {\rm e}^{h} \left( \left|\nabla^{\perp}X\right|^2 + \frac{h'(r)}{r}|X|^2 -  \sum_{i, j=1}^k \langle \textbf{B}(e_i,e_j),X \rangle ^2 \right) d \mathcal{H}_k \geq 0,
\end{equation}
since $\mathcal{C}_2$ has a globally flat normal bundle, there exists global normal parallel frame $e_{\alpha}(k+1\leq \alpha \leq n)$ on $\mathcal{C}_2$, we set $X=\sum\limits_{\alpha=k+1}^{n} \xi^{\alpha}e_{\alpha}$, then $\left|\nabla^{\perp}X\right|^2= \sum\limits_{\alpha=k+1}^{n} \left| \nabla \xi^{\alpha} \right|^2 $ and $|X|^2= \sum\limits_{\alpha=k+1}^{n} \left( \xi^{\alpha} \right)^2$, note that 
$$
\sum_{i, j=1}^k \langle \textbf{B}(e_i,e_j),X \rangle ^2= \sum_{i, j=1}^k \left( \sum_{\alpha=k+1}^{n} \xi^{\alpha} \langle \textbf{B}(e_i,e_j), e_{\alpha} \rangle  \right)^2 \leq  \sum_{\alpha=k+1}^{n} \left( \xi^{\alpha} \right)^2 \cdot |\textbf{A}|^2,
$$
we have 
$$
\delta^2 \mathcal{E}_f \geq \int_{\mathcal{C}_2} {\rm e}^{h} \sum_{\alpha=k+1}^{n} \left( \left| \nabla \xi^{\alpha} \right|^2  - |\textbf{A}|^2 \left( \xi^{\alpha} \right)^2 + \frac{h'}{|x|} \left( \xi^{\alpha} \right)^2 \right) d \mathcal{H}_{k},
$$
where $|\textbf{A}|^2$ denotes the square norm of the second fundamental forms of $\mathcal{C}$ in $\mathbb{R}^n$. Regardless of the normal indices, this is just the stability inequality for an oriented minimal hypercone. We want to prove that, under some suitable condition for $h$,
\begin{equation}\label{cmz16}
\int_{\mathcal{C}_2} {\rm e}^{h} \left( \left| \nabla \xi \right|^2  - |\textbf{A}(x)|^2 \xi^2 + \frac{h'}{|x|} \xi^2 \right) d \mathcal{H}_{k} \geq 0
\end{equation}
always holds for any $C^1$ function $\xi$ compactly supported in $\mathcal{C}_2$, then $\delta^2 \mathcal{E}_f \geq 0$, we can conclude $\mathcal{C}$ is $f$-stable.

\medskip\noindent

The stationary equation for $\mathcal{C}_2$ reads:
$$
\int_{\mathcal{C}_2} {\rm e}^{h} \left\{\operatorname{div} X+ \frac{h'}{|x|}\langle x, X\rangle  \right\} d \mathcal{H}_{k} =0, 
$$
for any $X\in C_{c}^1(\mathcal{C}_2,\mathbb{R}^n)$. Letting  $\xi \in C_c^1(\mathcal{C}_2,\mathbb{R})$ be arbitrary and put $X(x):=\frac{x}{|x|^2} \xi^2$. A direct calculation yields $\operatorname{div} X=\frac{(k-2)}{|x|^2} \xi^2+\frac{2 \xi}{|x|^2}(x \cdot \nabla \xi)$, then we have
$$
\int_{\mathcal{C}_2} {\rm e}^{h} \left\{(k-2+|x|h')|x|^{-2} \xi^2+2|x|^{-2} \xi(x \cdot \nabla \xi)\right\} d \mathcal{H}_{k}=0.
$$

Therefore we have
$$
\begin{aligned}
\int_{\mathcal{C}_2} \frac{k-2+|x|h'}{2} |x|^{-2} \xi^2 {\rm e}^{h} d \mathcal{H}_{k} & =-\int_{\mathcal{C}_2}  |x|^{-2} \xi(x \cdot \nabla \xi) {\rm e}^{h} d \mathcal{H}_{k} \\
& \leq\left\{\int_{\mathcal{C}_2}  |x|^{-2} \xi^2 {\rm e}^{h} d \mathcal{H}_{k}\right\}^{\frac{1}{2}}\left\{\int_{\mathcal{C}_2}  {\rm e}^{h} |\nabla \xi|^2  d \mathcal{H}_{k}\right\}^{\frac{1}{2}}
\end{aligned}
$$

In case the condition $th'(t)$ has a low bound $b >2-k$ on $\mathbb{R}^+$ is satisfied, we have
$$
\int_{\mathcal{C}_2}  {\rm e}^{h} |\nabla \xi|^2  d \mathcal{H}_{k} \geq \frac{k-2+b}{2} \int_{\mathcal{C}_2} \frac{k-2+|x|h'}{2} |x|^{-2} \xi^2 {\rm e}^{h} d \mathcal{H}_{k},
$$
thus if 
\begin{equation}\label{st1}
\begin{aligned}
\frac{(k-2+b)(k-2+|x|h'(|x|))}{4} & \geq |x|^{2} \left( |\textbf{A}(x)|^2-\frac{h'}{|x|} \right) \\
&=|\textbf{A}(\frac{x}{|x|})|^2-|x|h'(|x|) \\
&=|\textbf{A}_{\Sigma}|^2-|x|h'(|x|),
\end{aligned}
\end{equation}
then we can conclude that \eqref{cmz16} is non-negative, where a standard fact $|x|^{2} |\textbf{A}(x)|^2= |\textbf{A}(\frac{x}{|x|})|^2$ for $\mathcal{C}$ is used (for example, see \cite[Section 2.3]{JC22}). A more simplified condition than \eqref{st1} is 
\begin{equation}\label{st2}
\begin{aligned}
|\textbf{A}_{\Sigma}|^2 \leq \frac{(k-2+b)^2}{4}+b,
\end{aligned}
\end{equation}
where $|\textbf{A}_{\Sigma}|^2$ is the square norm of the second fundamental forms of $\Sigma \subset S^{n-1}(1)$. $\Box$

\medskip\noindent

Studying the Bernstein problem of high codimensional minimal cones is hard. Also for minimal hypercones with a general positive radial integral $f(x)={\ e}^{h(|x|)}$, the problem is not easy: when the standard arguments (see \cite[Appendix B]{Simon83} or the proof of \cite[Theorem 2.3]{DH23}) are  applied, the key equation \cite[Equation (15)]{DH23} will be written as:
$$
2\int_{\mathcal{C}_2} {\rm e}^{h}\xi^2 |\textbf{A}(x)|^2 \left(\frac{1-|x|h'}{|x|^2} \right) d \mathcal{H}_{n-1} \leq 
\int_{\mathcal{C}_2} {\rm e}^{h} |\textbf{A}(x)|^2  \left| \nabla \xi \right|^2 d \mathcal{H}_{n-1}.
$$
Also, restricted that $\xi=\xi(r)(r=|x|)$, it leads:
$$
\int_{\Sigma} |\textbf{A}_{\Sigma}|^2 d \mathcal{H}_{n-2} \cdot  \left(\int_{0}^{\infty} \left(\frac{d\xi}{dt} \right)^2 r^{n-4} {\rm e}^{h}dr -  \int_{0}^{\infty} \xi^2 (1-rh'(r))r^{n-6} {\rm e}^{h}dr \right) \geq 0,
$$
when we choose $\xi(r)=r^{\alpha}$ for $r \leq 1$ and $\xi(r)=r^{\beta}$ for $r \geq 1$, to conclude under what condition for $\alpha,\beta$ such that
\begin{equation}\label{te}
    \int_{0}^{\infty} \xi^2 (1-rh'(r))r^{n-6} {\rm e}^{h}dr=\int_{0}^{1} (1-rh'(r))r^{2\alpha+n-6} {\rm e}^{h(r)}dr + \int_{1}^{\infty} (1-rh'(r))r^{2\beta+n-6} {\rm e}^{h(r)}dr <\infty,
\end{equation}
we find that for general $h(r)$, from \eqref{te}, we can not conclude some concrete valid conditions such that \eqref{te} is satisfied. Then, the study for stable Bernstein problems is omitted here.

\medskip\noindent

We now give a condition for $f(x)={\rm e}^{h(|x|)}$, such that a stable minimal hypercone must be $f$-stable, following the arguments of Ding \cite[Section 3.3]{D20}.

\medskip\noindent

\textbf{Theorem 2.11 (Theorem 1.8 in Introduction)}. \textit{If for any $t>0$, the ordinary differential inequality $t(h')^2+2th''+2nh'\geq 0$ is satisfied for the radial integrals: $f(x)={\rm e}^{h(|x|)}$, then a stable, oriented minimal hypercone in $\mathbb{R}^n$ is also $f$-stable.} 

\medskip\noindent

\textbf{Proof}. Denote the compactly supported normal variation field on $\mathcal{C}_2$ by $X=\xi v$, $v$ is a global normal vector field of $\mathcal{C}_2$. Recall that
\begin{equation}
\delta^2 \mathcal{E}_f=
\int_{\mathcal{C}_2} {\rm e}^{h} \left( \left| \nabla \xi \right|^2  - |\textbf{A}(x)|^2 \xi^2 + \frac{h'}{|x|} \xi^2 \right) d \mathcal{H}_{n-1},
\end{equation}
by setting drifted Laplacian $\Delta_h:=\Delta+ \langle \nabla h, \nabla \cdot \rangle= {\rm e}^{-h} {\rm div} \left( {\rm e}^{h} \nabla \cdot \right)$, since $\xi$ is a compact supported function, following the divergence theorem, we can get 
\begin{equation}
\delta^2 \mathcal{E}_f=
\int_{\mathcal{C}_2} {\rm e}^{h} \xi \left( - \Delta \xi - \langle \nabla h, \nabla \xi \rangle - |\textbf{A}(x)|^2 \xi + \frac{h'}{|x|} \xi \right) d \mathcal{H}_{n-1}.
\end{equation}

We also let $\Delta_{\Sigma}$ be the Laplacian in $\Sigma$, for $\xi=\xi(y,r) \in C_{c}^{2}(\Sigma \times (0,\infty))$, the following relation is well-known (\cite{Sim68}):
$$
\Delta \xi =r^{-2} \Delta_{\Sigma} \xi + (n-2) r^{-1} \frac{\partial \xi}{\partial r} +\frac{\partial^2 \xi}{\partial r^2},
$$
combined with $\langle x, \nabla \xi \rangle = \langle x, D \xi \rangle = x( \xi)= r \frac{\partial}{\partial r}(\xi)$, we have 
$$
\delta^2 \mathcal{E}_f=
\int_{\Sigma \times \mathbb{R}^{+}}  \left( - \Delta_{\Sigma} \xi - 
\left|A_{\Sigma}\right|^2 \xi -(n-2) r \frac{\partial \xi}{\partial r} - r^2 \frac{\partial^2 \xi}{\partial r^2} -r^2 h'(r) \frac{\partial \xi}{\partial r} + rh'(r) \xi \right) \xi {\rm e}^{h} r^{n-4} d \mu_{\Sigma} d r, 
$$
since $\int_{\Sigma \times \mathbb{R}^{+}} \frac{\partial}{\partial r}(r^{n-2} \xi {\rm e}^{h} \frac{\partial \xi}{\partial r} )=0$, it reduce to 
\begin{equation}\label{svco}
\delta^2 \mathcal{E}_f=
\int_{\Sigma \times \mathbb{R}^{+}}  \left( - \xi \Delta_{\Sigma} \xi - 
\left|A_{\Sigma}\right|^2 \xi^2 + r^2  \left( \frac{\partial \xi}{\partial r} \right)^2+ rh'(r) \xi^2 \right) {\rm e}^{h} r^{n-4} d \mu_{\Sigma} d r. 
\end{equation}

Since
$$
\begin{aligned}
& \int_0^{\infty}\left(\frac{\partial}{\partial r}\left(\xi e^{\frac{h}{2}}\right)\right)^2 r^{n-2} dr= \int_0^{\infty}\left(\left(\frac{\partial \xi}{\partial r}\right)^2+h'(r)\xi \frac{\partial \xi}{\partial r}+\frac{(h')^2}{4} \xi^2 \right) r^{n-2} e^{h} d r \\
& \quad=\int_0^{\infty}\left(\frac{\partial \xi}{\partial r}\right)^2 r^{n-2} e^{h} d r- \frac{1}{4} \int_0^{\infty} \left( r(h')^2+2rh''+(2n-4)h' \right) \xi^2 r^{n-3} {\rm e}^{h} d r
\end{aligned}
$$

then finally \eqref{svco} writes as:

\begin{equation}\label{svco1}
\begin{aligned}
    \delta^2 \mathcal{E}_f &=
\int_{\Sigma \times \mathbb{R}^{+}}  \left( - \xi {\rm e}^{\frac{h}{2}} \Delta_{\Sigma} \left(\xi {\rm e}^{\frac{h}{2}}\right) - 
\left|A_{\Sigma}\right|^2 \left(\xi {\rm e}^{\frac{h}{2}}\right)^2 + r^2  \left( \frac{\partial}{\partial r} \left(\xi {\rm e}^{\frac{h}{2}} \right) \right)^2 \right) r^{n-4} d \mu_{\Sigma} d r \\
& \quad + \frac{1}{4} \int_{\Sigma \times \mathbb{R}^{+}} \left( r(h')^2+2rh''+2nh' \right) \xi^2 r^{n-3} {\rm e}^{h} d \mu_{\Sigma} d r
\end{aligned}
\end{equation}

Note the first integral term in \eqref{svco1} is just the stability inequality of area functional for $\xi {\rm e}^{\frac{h}{2}}$ (cf. \cite{Sim68}), hence, we have proved the result. $\Box$

\medskip\noindent

\textbf{Remark 2.12}. \textit{It is clearly that, if $t(h')^2+2th''+2nh'\leq 0$ on $\mathbb{R}^+$, then a $f$-stable minimal cone must be a stable minimal cone for the area-functional.}

\medskip\noindent
\medskip\noindent

\section{Lawlor's methods}
\subsection{Lawlor curvature criterion}

\medskip\noindent 

\textbf{A summary:} A regular minimal cone $\mathcal{C}$ satisfies Lawlor curvature criterion if there exists a well-defined (the normal wedges are disjoint, cf. \eqref{Normal radius}) area non-increasing  \textit{Lawlor retraction map} $\Pi$ from Euclidean space to $\mathcal{C}$ (see Definition 3.2) within a conical neighborhood of $\mathcal{C}$, the existence and area non-increasing properties of $\Pi$ ($J_{k} \Pi \leq 1$) both reduce to an ordinary differential inequality (written in polar coordinate $r=r(\theta)$) and the existence of some solutions of the associated O.D.E. tends to infinity: $r \rightarrow \infty$ as $\theta \rightarrow \theta_{0}$ (cf. Proposition 3.3 and below), the angle $\theta_0$ is called \textit{vanishing angle}. However, the non-intersection of $\theta_0$-normal wedges is not easy to check directly. A stronger condition, which is easy to calculate, states that if pointwisely two times of vanishing angle is no more than \textit{normal radius}, then $\mathcal{C}$ satisfy curvature criterion(cf. \cite[Theorem 1.2.1, Theorem 1.3.5]{La91}), hence be area-minimizing, and also weighted area-minimizing for any non-negative, non-decreasing radial integrals by Theorem 1.10.

\medskip\noindent 

For a given $k$-dimensional regular minimal cone $\mathcal{C}=C(\Sigma)$ in $\mathbb{R}^{n}$ with isolated singularity at the origin, we assume ${\rm dim} \ \Sigma =k-1\geq 2$, note that $\Sigma \subset S^{n-1}(1)$ is a smooth, minimal submanifold, and it is well-known that $\mathcal{C}$ must be flat when $k-1=1$. 

\medskip\noindent
  
\textbf{Definition 3.1. (\cite{La91})} \textit{\emph{(1)} For $p\in \Sigma$ and $\alpha \geq 0$, let $U_{p}(\alpha)$ be the union of all open normal geodesic (in the sphere) of length $\alpha$ from $p$. The \textit{normal wedge} $W_{p}(\alpha)$ is defined as the cone over $U_{p}(\alpha)$.}

\textit{\emph{(2)} Let $N_{p}:= {\rm max} \{\alpha | W_{p}(\alpha)\cap \mathcal{C}= \overrightarrow{op}\}$ , where $\overrightarrow{op} \in  \mathcal{C}$ is the ray through $p\in \Sigma$, we call $N_{p}$ the \textit{normal radius} at $p$.}

\medskip\noindent

Lawlor constructed a retraction map from some angle neighborhood of $\mathcal{C}\subset \mathbb{R}^{n}$ to the cone $\mathcal{C}$: he choose a point $p\in \Sigma$ and a unit normal vector $v$ of $\mathcal{C}$, select a curve in polar coordinate $(r(\theta),\theta)$ on the $\mathbf{2}$-dimensional subspace spanned by $\{\overrightarrow{op}, v\}$ of $\mathbb{R}^{n}$:
\begin{equation}\notag
 S_{p}=r(\theta) ( \cos  \theta \ \overrightarrow{op}+ \sin  \theta \ v),
\end{equation}
then extend $S_{p}$ homothetically on $\mathcal{C}$. 

\medskip\noindent

It is easy to check that Lawlor's construction is equivalent to giving a map from $\mathcal{C}-\{0\}$ to $\mathbb{R}^n$ for angle $\theta$-slices (whether $\theta$ depends on the moving of $x$ will not affect our results):
\begin{equation}\label{cuirewrite}
\Phi(x)=r(\theta)\left(\cos \theta \ x+\sin \theta  \ |x|  \widetilde{V}\left(\frac{x}{|x|}\right)\right),
\end{equation}
where $\widetilde{V}$ is an extended normal vector field of $v$ along the unit sphere.

\medskip\noindent
 
Lawlor's idea is that trying to find some function $r(\theta)$ such that both the $k$-dimensional Jacobian of $\Phi$ is always no less than $1$ thus area non-decreasing, and $\Phi$ extends to infinity when $\theta$ approach to some fixed angle $\theta_{0}$ (the following solutions of vanishing angle type, the angle is less than the focal radii, see \cite{PT88}), i.e. $\Phi_{p}$ takes values at a normal wedge of $\mathcal{C}$ at $p$.
And, if additionally, pointwisely the normal wedges $ W_{p}(\theta_{0}(p))$ disjoint: for any $p,q \in \Sigma, p\neq q$,
\begin{equation}\label{Normal radius}
 W_{p}(\theta_{0}(p))\cap W_{q}(\theta_{0}(q))=\varnothing, 
\end{equation}
then we can get a inverse map $\Pi:=\Phi^{-1}$ which is well-defined in an angle neighborhood denoted by $W:= \bigcup_{p\in \Sigma} W_{p}(\theta_{0}(p))$, define that $\Pi$ maps every point outside $W$ to the origin. Easy to see $\Pi$ is area non-increasing. Also, from \eqref{cuirewrite}, we can see that $\Pi$ is just the nearest-point projection up to a scalar $\frac{1}{r(\theta) \cos \theta}$.

\medskip\noindent
 
\textbf{Definition 3.2.} \textit{We call the above well-defined area non-increasing retraction map: $\Pi$ in an angle neighborhood of a given cone $\mathcal{C}$, the \textit{Lawlor retraction map}.}

\medskip\noindent

We now talk about the existence of Lawlor retraction map, in \cite[p.10-12]{La91}, Lawlor proved that\footnote{In \eqref{cuirewrite}, by varying $x$ in different tangential directions, readers can directly estimate the $k$-dimensional Jacobian of $\Phi$, which is the reciprocal of $J_{k}\Pi$.}:

\medskip\noindent

\textbf{Proposition 3.3. (the $k$-dimensional Jacobian of $\Pi$)}: \textit{the $k$-dimensional Jacobian of $\Pi$ on $\mathbb{R}^n$: $J_{k}\Pi \leq 1$ is equivalent to the following ordinary differential inequality 
\begin{equation}\label{va1}
    \frac{dr}{d\theta}\leq r \sqrt{r^{2k}({\rm cos}\theta)^{2k-2}\left(\mathrm{inf}_{v}{\rm det}(I-{\rm tan} \theta A_v)\right)^{2}-1}
\end{equation}
for any normal direction $v$, where $A_v$ is the shape operator of $\Sigma \subset S^{n-1}$ in normal direction $v$.}

\medskip\noindent

As pointed out by Lawlor, the method of curvature criterion is firstly discovered in the context of \textit{calibrations}, rather than retractions. And the area-nonincreasing requirements and the associated ordinary differential inequality above turn out to be equivalent to the comass of a \textit{vanishing calibration} is no larger than one, see \cite[Section 2]{La91}.

\medskip\noindent

To make the following vanishing angle $\theta_0$ as small as possible (the possibility for the non-intersection of $\theta_0$-normal wedges is bigger),  we replace the inequality in \eqref{va1} with equality. The solution of the corresponding ODE is analyzed by Lawlor: 
\begin{equation}\label{Vanishing angle}
\begin{cases}
\frac{dr}{d\theta}= r \sqrt{r^{2k}({\rm cos}\theta)^{2k-2}\left(\mathrm{inf}_{v}{\rm det}(I-{\rm tan} \theta A_v)\right)^{2}-1}\\
r(0)=1,
\end{cases}
\end{equation}
either $\frac{dr}{d\theta}$ vanish at some positive $\theta(p)$ or $r\rightarrow \infty$ as $\theta \rightarrow \theta_{0}(p)$, in the latter case, Lawlor can construct his retraction map $\Pi$ if in addition \eqref{Normal radius} is satisfied for the normal wedges with the angle $\theta_{0}$, then it follows that the cone is area-minimizing, that is the origin version of Lawlor curvature criterion (cf. \cite[Theorem 1.2.1]{La91}). And, we say a minimal cone \textit{satisfying Lawlor curvature criterion} (hence area-minimizing) if it satisfies the above steps. It is worth mentioning that for zero curvature cases (planes through the origin), there also exists the Lawlor retraction maps, see first line in Lawlor's table (cf. \cite[p.20]{La91})\footnote{A strict proof can be given by Proposition 3.5, readers should find the estimated vanishing angles and check that the normal radius of planes through the origin is equal to $\pi$.}.

 \medskip\noindent
 
\textbf{Definition 3.4.([La91])} \textit{We call the smallest $\theta_{0}(p)$ among all the normal directions $v$ the vanishing angle at $p$.}

\medskip\noindent

Lawlor reduced the condition \eqref{Normal radius} to the computation of normal radius, and he gave the following simplified---a weaker version of Lawlor curvature criterion that is more frequently used.

\medskip\noindent

\textbf{Proposition 3.5.(\cite{La91})} \textit{If the vanishing angle $\theta_{0}(p)$ exists for any points $p\in \Sigma$ and pointwisely satisfy
\begin{equation}\label{two times vanishing angle}
2\theta_{0}(p)\leq N(p),
\end{equation}
where $N(p)$ is the normal radius at $p$, then the cone $\mathcal{C}=C(\Sigma)$ is area-minimizing (in the sense of mod $2$ when $\Sigma$ is nonorientable).}

\medskip\noindent
    
\textbf{Remark 3.6.} \textit{Lawlor had computed the estimated vanishing angles---the approximate solution of \eqref{Vanishing angle}, for $dim \ \mathcal{C} \leq 12$ which depending only on $sup_{v}||A_{v}||^2$ and $dim \ \mathcal{C}$, he also gave an upper bound formula of vanishing angle for the cases of $dim \ \mathcal{C} \geq 13$. These works make Curvature Criterion very convenient to use, see \cite{Ke94},\cite{Ka02},\cite{XYZ18},\cite{TZ20},\cite{JC22},\cite{JCX22}, etc.}

\medskip\noindent

In \cite[Sections 3,4]{La91}, Lawlor conducted more in-depth discussions for his curvature criterion. For plenty of minimal cones that admit global parallel normal vectors (such as with globally flat normal bundles), the curvature criterion is both sufficient and necessary. This means that if there are no solutions of vanishing angle type (for some special class of geometric cones, the $\theta_{0}$ normal wedges are naturally disjoint), then we can find a comparison surface with less area than the cone, and the cone is not area-minimizing. Such cones include the minimal isoparametric hypercones (including all homogeneous minimal hypercones), the cones over products of spheres, and the cones over principle orbits of polar group actions. Based on these arguments, Lawlor was able to prove that: the cone $\mathcal{C}_{1,5}$ (proved one-sided area-minimizing by Lin \cite{Lin87}) and the cone over $(SO(2) \times SO(8))/ (\mathbb{Z}_2 \times SO(6))$ are stable but not area-minimizing.

\medskip\noindent
\medskip\noindent 

\subsection{Lawlor directed slicing}

\medskip\noindent 

For proof of our non-regular $f$-minimizing cones, we introduce Lawlor's directed slicing method. The method was set up in \cite{La98} by Lawlor, where comprehensive investigations and research are given. Roughly speaking, it turns an optimization problem into a family of localized optimization problems. We quote his original words:

 \textit{Directed slicing is particularly advantageous for proving area-minimization of union of minimal surfaces, and of highly singular or unorientable surfaces, as well as slight perturbations of such surfaces. Questions of this type have previously been very challenging.}

As its applications, two significant families of minimal real matrix varieties: determinantal varieties and Pfaffian varieties are proved to be area-minimizing stratified cones except for several cases (\cite{KL99},\cite{CJX24}). 

\medskip\noindent 

{\bf Definition 3.7} (cf. \cite[Definition 2.1]{KL99})\,\, \textit{Let $M$ be a $k$-dimensional surface (rectifiable current or rectifiable current reduced modulo 2) in $\mathbb{R}^{n}$. For $0\leq d<k$, a \textit{$d$-dimensional slicing} of $M$ is a collection of pairwise disjoint, $d$-dimensional rectifiable subsets of $M$, if the union of these rectifiable subsets covers $M$, $\mathcal{H}^{k}$ a.e., then the slicing is called \textit{full}. 
We are concerned with those $(n-k+d)$-dimensional slicings of $\mathbb{R}^{n}$, which its intersections with $M$ consist of $d$-dimensional slicings of $M$, then we call those $d$-dimensional subsets of $M$: \textit{slices}, and those $(n-k+d)$-dimensional subsets of $\mathbb{R}^{n}$: \textit{slicing sets} or just \textit{slicings}.}

\medskip\noindent

Those slicing sets of the ambient space $\mathbb{R}^{n}$ can be extended from the slices of $M$, i.e. the normal spaces slicings, see \cite[Section 3]{KL99}. Such an extended slicing is valid locally: the slicing sets are pairwise disjoint in a neighborhood of $M$ bounded by the focal sets (cf. \cite{PT88}).

\medskip\noindent 

Often, the slicing sets of $\mathbb{R}^{n}$ are given by level sets of some $C^1$ map: $f:\mathbb{R}^{n} \rightarrow \mathbb{R}^{k-d}$, the important \textit{weighting function} is defined to be the reciprocal of Jacobians (see \cite[4.6]{La98} and \cite[3.2.1]{Fe69}, or \cite[3.6]{M16}):
$$
w(x)=\frac{1}{(J_{k-d}f)|_{M}}, \ \ \ \  J_{k-d}f={\rm max}||(Df)_{p}(v_{1}\wedge\cdots\wedge v_{k-d})||,
$$
where $v_{1},\ldots,v_{k-d}$ are unit, orthogonal vectors of $\mathbb{R}^{n}$ at $p$, the restriction of $f$ on $M$ gives a new Jocabian $J_{k-d}(f|_{M})$ and clearly they satisfy  
$$
J_{k-d}(f|_{M})\leq (J_{k-d})f|_{M},
$$
the equality is attained when every level set of $f$ intersects $M$ orthogonally.

\medskip\noindent

Let $g$ be an integrable function on $M$ with the induced Hausdorff measure, and the area-coarea formula in GMT (see \cite[3.2.22]{Fe69} and \cite[Theorem 4.5.1]{La98}) states that
\begin{equation}\label{areacoarea}
\int_{M}g J_{k-d}(f|_{M})d \mathcal{H}^{k}=\int_{\mathbb{R}^{k-d}}\left(\int_{f^{-1}(y)\bigcap M}g\; d\mathcal{H}^{d}\right) d \mathcal{H}^{k-d}.
\end{equation}

Instead $g$ by the weighting function, then
$$
{\rm Area} M \geq \int_{M} \frac{1}{(J_{k-d}f)|_{M}} J_{k-d}(f|_{M})d \mathcal{H}^{k}
=\int_{\mathbb{R}^{k-d}}\left(\int_{f^{-1}(y)\bigcap M}\frac{1}{(J_{k-d}f)|_{M}} d\mathcal{H}^{d}\right) d \mathcal{H}^{k-d}.
$$

\medskip\noindent

\textbf{Remark 3.8}. \textit{Follow \cite[3.2.22]{Fe69}, every slice is rectifiable.}

\medskip\noindent

A key point in Lawlor's treatise is that he required the slicing sets to intersect $M$ orthogonally, then the above inequality turns out to be equality, and the inside integrals suggest a new weighted area-minimizing question on each slicing set, if we slice $\mathbb{R}^{n}$ twice or more, then the inside weighting functions will be the product of the weighting functions produced by iterated slicing every time (see \cite[4.5.3]{La98}). 

\medskip\noindent

Every composition function $h=g \circ f$ gives the same slicing defined by the level sets, and they are different up a constant on each level set, which will not affect the integral result (see \cite[Proposition 4.6.4]{La98} and its proof).

\medskip\noindent

\textbf{Remark 3.9.} \textit{In the context of Directed Slicing, Lawlor retraction map or its inverse map: $\Phi$ in \eqref{cuirewrite}, is $(n-k)$-dimensional slicings over the $0$-dimensional slices (points) of $\mathcal{C}$, and the fibers are also orthogonal to $\mathcal{C}$ (in \eqref{Vanishing angle}, $\frac{dr}{d\theta}(0)=0$). So, $J_{k}\Pi \leq 1$ just means that the weighting function attains minimum on the cone, the associated O.D.in.E. leads to the computation of vanishing angle in Curvature Criterion. We note that, for the applications of Directed Slicing, one also needs to ensure that slicing sets are disjoint (for example, see \cite[Proposition 6.3]{KL99} and \cite[Proposition 4.2]{CJX24}), and for regular cones, it is just the condition \eqref{Normal radius}.}

\medskip\noindent
\medskip\noindent

\section{$f$-minimizing cones}
\subsection{$f$-minimizing regular cones and $f$-minimizing planes}

\medskip\noindent 

In this subsection, we prove Theorem 1.10. For regular cones, more precisely, we want to prove a regular minimal cone $\mathcal{C}$ satisfying Lawlor curvature criterion is $f$-minimizing, among rectifiable currents having the same boundary with $\mathcal{C}_1$ or among rectifiable currents reduced modulo $2$ having same mod-2 boundary with $\mathcal{C}_1$, for the non-negative radial integrals: $f(x)=g(|x|)$ such that $g$ is non-decreasing. 

\medskip\noindent 

\textbf{Proof:} \textbf{I. $f$-minimizing regular cones:} As before, let $\mathcal{C}_{1}$ denote the part of $\mathcal{C}$ inside the unit ball, $S$ a comparision surface. Since $\mathcal{C}$ satisfies Lawlor curvature criterion, there exist normal wedges, denoted by $W$, i.e. an angle neighborhood of the cone $\mathcal{C}$:
$$
W= \bigcup_{p\in \Sigma} W_{p}(\theta_{0}(p)) \subset \mathbb{R}^{n},
$$
where $W_{p}(\theta_{0}(p))$ are disjoint sets for different points $p$ in the link of $\mathcal{C}$. Within $W$---the open subset of $\mathbb{R}^{n}$, there exists a Lawlor retraction map $\Pi$---the nearest point projection up to a scalar $\frac{1}{r(\theta) \rm{cos} \theta}$, $\Pi:W\rightarrow \mathcal{C}_{1}$.

\medskip\noindent

By general area-coarea formula \eqref{areacoarea}, we have 
\begin{equation}\label{beq}
\begin{aligned}
\int_{S \cap W }f(x) J_{k}\Pi|_{S \cap W}(x) d \mathcal{H}_k (x)  &= \int_{\mathcal{C}_1}  \left( \int _{\Pi^{-1}(y) \cap (S \cap W) } f(x) d\mathcal{H}_0 \right)d \mathcal{H}_k (y) \\
&= \int_{\mathcal{C}_1}\left( \sum_{x \in \Pi^{-1}(y) \cap S} f(x) \right) d \mathcal{H}_k(y)   \\
&\geq \int_{\mathcal{C}_{1}} f(y) d \mathcal{H}_k(y)=\mathcal{E}_f(\mathcal{C}_{1})
\end{aligned}
\end{equation}
the second equality follows from the fact that the inverse images of $\Pi$ are contained in $W$. To prove the inequality, we first need to show almost every slicing $\Pi^{-1}(y)$ intersects $S$, this fact had already been proven in the proofs of \cite[Theorem 1.2.1]{La91} and \cite[Theorem 8.3]{KL99}. The former is based on a variant of the Constance Theorem in GMT, we note that in the second one, the cases of currents reduced mod 2 were omitted, and we discuss them here: 

\medskip\noindent

(1)~~(rectifiable currents): we follow the proof in \cite[Theorem 8.3]{KL99}: if $\partial S= \partial \mathcal{C}_1$, then $S$ and $\mathcal{C}_1$ bounds a $(k+1)$-dimensional current, denoted by $R$, i.e., $\partial R=S-\mathcal{C}_{1}$, since the slicing commutes with the boundary operator, up to sign (see 4.3 in \cite[4.3]{Fe69}), it follows that
$$
(S-\mathcal{C}_{1})\cap \Pi^{-1}(y) = (\partial R) \cap \Pi^{-1}(y) = \pm \partial (R \cap \Pi^{-1}(y)),
$$
for $\mathcal{H}_k$ a.e. $y\in \mathcal{C}_1$, $R \cap \Pi^{-1}(y)$ is a curve or union of curves, and its boundary will consist of at least two points, which implies that the sets $\Pi^{-1}(y) \cap S$ are non-empty.

\medskip\noindent 

(2)~~(rectifiable currents reduced mod 2): If $\partial S =\partial \mathcal{C}_1 \ (mod \ 2)$, then there exist currents $R, Q$, such that $S-\mathcal{C}_{1}=\partial R+2Q$, it follows that
$$
(S-\mathcal{C}_{1})\cap \Pi^{-1}(y) = (\partial R) \cap \Pi^{-1}(y) + 2 Q \cap \Pi^{-1}(y) = \pm \partial (R \cap \Pi^{-1}(y)) +2 Q \cap \Pi^{-1}(y),
$$
for $\mathcal{H}_k$ a.e. $y\in \mathcal{C}_1$, $R \cap \Pi^{-1}(y)$ is a curve or union of curves, then $\partial (R \cap \Pi^{-1}(y))$ is a zero-dimensional current and a boundary---that is, an even number of points (counting multiplicity) (cf. \cite[Theorem 10.1]{M16}). Thus, the number of points on the right-hand side is even, and $\mathcal{C_1} \cap \Pi^{-1}(y)=y$ consists of a single point which also implies that almost every slicing $\Pi^{-1}(y)$ intersect $S$ in the cases of rectifiable currents mod 2.

\medskip\noindent 

The inequality then is proved by our assumption: $f(x)=g(|x|)$ is radially non-decreasing and clearly $|\Pi^{-1}(y)|=|\Phi(y)|=r|y|>|y|$ for every $y \in \mathcal{C}_1$. 

\medskip\noindent

Note that $J_{k}\Pi|_{S \cap W}(x) \leq J_{k}\Pi(x)$ for every $x \in S \cap W$, and also the construction of Lawlor retractions yields that (see Proposition 3.3): 
$J_{k}\Pi(x) \leq 1$. 

\medskip\noindent

Since $f$ is non-negative on $\mathbb{R}^{n}$, then we can conclude that 
$$
\mathcal{E}_f(\mathcal{C}_{1}) \leq \mathcal{E}_f(S \cap W) \leq \mathcal{E}_f(S),
$$
for every comparison surface with the same boundary with $\mathcal{C}_1$. 

\medskip\noindent 

\textbf{II. $f$-minimizing planes:} for the proof of $\mathbb{R}^{k} \subset \mathbb{R}^{n}$ minimizing $\mathcal{E}_ f$ for any radially non-decreasing integrals (seen as flat cones), firstly, for the cases of $k\geq 3$, there also exists the angle neighborhoods and Lawlor retractions described as before, and they satisfy curvature criterion, see Chapter 1 in \cite{La91}.

\medskip\noindent 

For a proof of all $\mathbb{R}^{k}$ and all possible domain $\Omega \subset \mathbb{R}^{k}$ (not only the part inside the unit ball), in \eqref{beq}, we can just choose the retraction map to be the orthogonal projection from $\mathbb{R}^{n}$ to $\mathbb{R}^{k}$ and replace $\mathcal{C}_1$ to be $\Omega$, then a direct application of area-coarea formula yields the proof, similar to \textbf{I}. $\Box$

\medskip\noindent 
\medskip\noindent 

\subsection{$f$-minimizing determinantal varieties and Pfaffian varieties}

\medskip\noindent 

\textbf{4.2.1. An overview on slicing processes for determinantal varieties and Pfaffian varieties:} For determinantal varieties and Pfaffian varieties, they both need two times of slicings: primary slicings and secondary slicings, see \cite[Sections 6,7]{KL99} and \cite[Sections 4,5]{CJX24}.

\medskip\noindent

1. \textit{Determinantal varieites} $C(p,q,r) \subset \mathbb{R}^{pq}$: let $X_{ij}$ denote the $p$ by $q$ matrix with $1$ in the $(i,j)$ position and zeros everywhere, denote
$$
\widetilde{E}_1(x_{1},\ldots,x_{r})=\sum_{i=1}^{r}x_{i}X_{ii} \in C(p,q,r),
$$
where $x_{1}> \cdots > x_{r}$ are nonascending positive numbers.

\medskip\noindent

Let $K_1$ be the subspace spanned by vectors $X_{ij}$ for which $i>r,j>r$, $K_1$ is the normal space of $\widetilde{E}_1$.

Let $H_1=\{M_1+N_1: M_1 \in \widetilde{E}_1, N_1 \in K_1, and \ |N_1|<x_{r}\}$, where $x_r$ denotes the minimum nonzero singular value of $M_1$, it can be seen as "adding" the normal spaces on $\widetilde{E}_1$ bounded by focal radii at every regular point.

Recall \cite[Definition 6.1]{KL99}, the \textit{primary slices} of $C(p,q,r)$ are defined as the distinct images of adjoint actions on $\widetilde{E}_1$: $P\widetilde{E}_1Q^{T}$, $P\in SO(p), Q \in SO(q)$. $C(p,q,r)$ can be seen as the union of rotation images of the closure of $\widetilde{E}_1$.

The \textit{primary slicing sets} (or \textit{primary slicings}) of $\mathbb{R}^{pq}$ are defined as the distinct images of adjoint actions on $H_1$: $PH_1Q^{T}$, $P\in SO(p), Q \in SO(q)$, and they are proved to be disjoint in \cite[Proposition 6.3]{KL99}.

\medskip\noindent

Kerckhove and Lawlor computed the primary weighting function at regular points belonging to $\widetilde{E}_1$ (see \cite[Proposition 6.7]{KL99}): let $M_0=M\left(x_1, \ldots, x_r\right)=M(\vec{x})$ be a point of $\widetilde{E}_1$, let $v_1=$ $\left[\begin{array}{ll}0 & 0 \\ 0 & B_1\end{array}\right]$, $|B_1|=1$, be a unit normal to $C(p, q, r)$ at $M(\vec{x})$. For $|t|<x_r$, the weighting function $w_1$ for the primary slicing of $\mathbb{R}^{p q}$ satisfies
\begin{equation}\label{primary 1}
    w_1(M(\vec{x})+t v_1) \geqslant \prod_{1 \leqslant i<j \leqslant r}\left(x_i^2-x_j^2\right) \prod_{i=1}^r x_i^{p+q-2 r-2} \prod_{i=1}^r\left(x_i^2-t^2\right)
\end{equation}
with equality if and only if $t=0$ or $rank \ B_1=1$.

\medskip\noindent

The factors $(x_{i}^2-t^2)$ in \eqref{primary 1} give the motivation to define the secondary slicings which are 1-dimensional, they will decompose the wedge-shaped primary slicing sets into disjoint hyperbolic curves (cf. Figure 2 in \cite{La98}): set
$$
h_i(M(\vec{x})+t v_1)= \frac{x_i^2-t^2}{2}, \quad i=1, \ldots, r,
$$
and let the secondary slicing of $H_1$ be the collection of level sets of the function $h$, denoted them by $H_{c}$,
$$
H_{c}:=h^{-1}(c),
$$
where $c=(c_{1}^2/2,\ldots,c_{r}^2/2)$. Moreover, the union of secondary slicings sets $H_{c}$ covers $H_1$ a.e., and every $H_{c}$ intersects orthogonally with $\tilde{E}_1$ at a single point. 

\medskip\noindent

Lawlor computed the composite weighting function (for the case $r=1$, the terms $(x_{i}^2-x_{j}^2)$ disappeared):
\begin{equation}\label{cf1}
 w_{1}w_{2}(M(\overrightarrow{x})+tv_1) \geq \frac{\prod\limits_{1 \leqslant i<j \leqslant r}\left(x_{i}^2-x_{j}^2\right)\prod\limits_{i=1}^{r}x_{i}^{p+q-2r-3}\prod\limits_{i=1}^{r}\left(x_{i}^2-t^2\right)}{\sqrt{1+\sum\limits_{i=1}^{r}\frac{t^2}{x_{i}^{2}}}},  
\end{equation}
with equality if and only if $t=0$ or $rank \ B_1=1$.

\medskip\noindent
\medskip\noindent

2. \textit{Pfaffian varieties} ${\bf C}(m,2r) \subset Skew(m,\mathbb{R})\cong \mathbb{R}^{\frac{m(m-1)}{2}}$: Let $Y_{ij}(i<j)$ denote the $m$ by $m$ skew-symmetric matrix with $1$ in the $(i,j)$ position, -$1$ in the $(j,i)$ position and zeros everywhere, denote
$$
\tilde{E}_2(y_{1},\ldots,y_{r})=\sum_{i=1}^{r}y_{i}Y_{2i-1,2i} \in {\bf C}(m,2r),
$$
where $y_{1}> \cdots > y_{k}$ are nonascending positive numbers.

\medskip\noindent

Let $K_2$ be the subspace spanned by vectors $Y_{ij}$ for which $2r<i<j\leq m$, $K_2$ is the normal space of ${\bf C}(m,2r)$ at regular points.

Let $H_2=\{M_2+N_2: M_2 \in \widetilde{E}_2,N_2 \in K_2, and \ |N_2|<y_{2r-1,2r}\}$, where $y_{2r-1,2r}$ denotes the minimum nonzero singular value of $M_2$.

Recall \cite[Definition 4.1]{CJX24}, the \textit{primary slices} of ${\bf C}(m,2r)$ are defined as the distinct images of adjoint actions on $\widetilde{E}_2$: $P\widetilde{E}_2P^{T}$, $P\in SO(m)$.
${\bf C}(m,2r)$ can be seen as the union of rotation images of the closure of $\widetilde{E}_2$.

The \textit{primary slicing sets} of $Skew(m,\mathbb{R})$ are defined as the distinct images of adjoint actions on $H_2$: $QH_2Q^{T}$, $Q \in SO(m)$, and they are also proved to be disjoint in \cite[Proposition 4.2]{CJX24}.

\medskip\noindent

In \cite{CJX24}, similar expressions for primary weighting function $w_1$ and the same definition for the secondary slicings are given, the composite weighting functions $w_1w_2$ is:

\begin{equation}\label{cf2}
 w_{1}w_{2}(M(\overrightarrow{y})+tv_2) \geq \frac{\prod\limits_{1\leq i <j \leq r}\left(y_{i}^2-y_{j}^2\right)^2\prod\limits_{i=1}^{r}y_{i}^{2m-4r-5}\prod\limits_{i=1}^{r}(y_{i}^2-t^2)}{\sqrt{1+\sum\limits_{i=1}^{r}\frac{t^2}{y_{i}^{2}}}}.  
\end{equation}
where $M(\overrightarrow{y})=\sum\limits_{i=1}^{r}y_{i}Y_{2i-1,2i}$ and $ y_1 > \cdots >y_{r}>|t|$, $v_2=\left[\begin{array}{ll}0 & 0 \\ 0 & B_2\end{array}\right]$ is a unit normal of ${\bf C}(m,2r)$, the equality is satisfied if and only if $t=0$ or $tr B_2 B_2^{T}=2$. In particular, the equality is satisfied for Pfaffian hypersurfaces ${\bf C}(2n,2n-2)$.

\medskip\noindent
\medskip\noindent 

\textbf{4.2.2. Compensating for the composite weighting functions:} A key step in the proof of area-minimization of these algebraic varieties (cones) (cf. \cite[Theorem 8.3]{KL99} and \cite[Proposition 6.1]{CJX24}) is that ones need to show the composite weighting functions $w_1w_2$ within each secondary slicing set attains its minimum value at these cones.

\medskip\noindent 

Thus, if we want to prove those minimal cones are also $f$-minimizing for the Euler-Dierkes-Huisken variational problem, a reasonable chosen condition for $f(x)=g(|x|)$ is that, along every secondary slicing $H_{c}$,  the new weighting functions: $fw_1w_2$ attains its minimum values at those cones, i.e. 
\begin{equation}\label{fcd}
fw_1w_2(H_{c} \cap C_{1}) \leq fw_1w_2(H_{c}).
\end{equation}
We should note that, besides the establishment of \eqref{fcd} here, many necessary lemmas and techniques had been given in: \cite{La91},\cite{La98},\cite{KL99}, etc.

\medskip\noindent

Recall that, the composite weighting function: $w_1w_2$ attains minimum values at the cones: $C(p,q,r)(p+q-2r\geq 4)$ and  ${\bf C}(m,2r)(m-2r\geq 3)$, thus for those area-minimizing cones proved in \cite{KL99} and \cite{CJX24}, we can just choose the radial integrals: $f(x)=g(|x|)$ such that $g(t)$ is non-decreasing. 

\medskip\noindent

The left determinantal varieties: $C(p,q,r)(p+q-2r=3)$ are just $C(p,p+1,p-1)$, $C(p,q,r)(p+q-2r=2)$ are just determinantal hypersurfaces $C(p,p,p-1)$. The left Pfaffian varieties are just Pfaffian hypersurfaces ${\bf C}(2n,2n-2)$. Whether they are area-minimizing is still unsolved! In the following, we will prove that for some of them, there do exist radial integrals such that the \eqref{fcd} is satisfied. However, for generic ones, we prove that there can \textbf{NOT} exist any $C^1$ radial integral (away from the origin) such that \eqref{fcd} is satisfied!

\medskip\noindent 


\medskip\noindent 

\textbf{Proposition 4.1}  \textit{\emph{(1)} For every non-negative radial integral function: $f(x)=g(r), r=|x|$, if $g$ is non-decreasing, then the minimum values of the new weighting function $fw_1w_2$ within every secondary slicing set equals to its values taken at those points belong to the area-minimizing cones $C(p,q,r)(p+q-2r\geq 4)$ and ${\bf C}(m,2r)(m-2r\geq 3)$;}

 \textit{\emph{(2)} For determinantal hypersurface $C(2,2,1)=C(S^1 \times S^1)$, determinantal variety $C(2,3,1)$ and Pfaffian hypersurface ${\bf C}(4,2)=C(S^2 \times S^2)$, given any radial integral functions such as $\frac{g(r)}{r}$ is non-decreasing, then within each secondary slicing set, $fw_1w_2$ takes its minimum at the cones, those integral functions include $f(x)=g(|x|)=|x|^{\alpha}(\alpha \geq 1)$ etc;} 

 \textit{\emph{(3)} For the left  determinantal varieties: $C(n,n+1,n-1)(n\geq 3)$, generic determinantal hypersurfaces $C(n,n,n-1)(n\geq 3)$ and the left generic Pfaffian hypersurfaces: ${\bf C}(2n,2n-2)(n\geq 3)$, there \textbf{doesn't} exist any $C^1$ radial integrals such that along all the secondary slicings, $fw_1w_2$ takes its minimum at the cones.}

\medskip\noindent 

{\bf Proof.} \,\, The proof of (1) is trivial. We aim to prove (2),(3).

\medskip\noindent

Recall \eqref{cf1}\footnote{the condition $rank \ B=1$ is naturally satisfied for $C(n,n+1,n-1)$ and $C(n,n,n-1)$!}, for determinantal varieties $C(n,l,n-1)$($l=n$ or $ l=n+1$), the composition weighting function $w_1w_2$ just equals to (for the case $n=2$, the terms $(x_{i}^2-x_{j}^2)$ disappeared): 
$$
w_{1}w_{2}(M(\overrightarrow{x})+tv_1)=\frac{\prod\limits_{1\leq i <j \leq n-1}(x_{i}^2-x_{j}^2)\prod\limits_{i=1}^{n-1}x_{i}^{l-n-1}\prod\limits_{i=1}^{n-1}(x_{i}^2-t^2)}{\sqrt{1+\sum\limits_{i=1}^{n-1}\frac{t^2}{x_{i}^{2}}}},
$$
where $M(\overrightarrow{x})=\sum\limits_{i=1}^{n-1}x_{i}X_{ii}, x_1 > \cdots >x_{n-1}>|t|$, and $v_1$ is a unit normal of $C(n,l,n-1)$.

\medskip\noindent

Recall \eqref{cf2}, for Pfaffian hypersurfaces ${\bf C}(2n,2n-2)$, the values of composition weighting function $w_1w_2$ equals (for $n=2$, the terms $(y_{i}^2-y_{j}^2)^2$ disappeared): 
$$
w_{1}w_{2}(M(\overrightarrow{y})+tv_2)=\frac{\prod\limits_{1\leq i <j \leq n-1}(y_{i}^2-y_{j}^2)^2\prod\limits_{i=1}^{n-1}(y_{i}^2-t^2)^2}{\prod\limits_{i=1}^{n-1}y_{i}\sqrt{1+\sum\limits_{i=1}^{n-1}\frac{t^2}{y_{i}^{2}}}},
$$
where $M(\overrightarrow{x})=\sum\limits_{i=1}^{n-1}y_{i}Y_{2i-1,2i}$ and $ y_1 > \cdots >y_{n-1}>|t|$, $v_2$ is a unit normal of ${\bf C}(2n,2n-2)$.

\medskip\noindent

Note that, the inner products are defined, and it is easy to check that the distances from the origin to those points $M(\overrightarrow{x})+tv_1$ equal to $\sqrt{x_{1}^2+ \cdots + x_{n-1}^2+t^2}$, to those points 
$M(\overrightarrow{y})+tv_2$ equals to $\sqrt{y_{1}^2+ \cdots+y_{n-1}^2+t^2}$. We will consider the questions for determinantal varieties and Pfaffian varieties uniformly, by identifying \textbf{$y_i$ with $x_i$}, etc. Then, along every fixed secondary slicing set $H_{c}(c=(c_{1}^2/2,\ldots,c_{n-1}^2/2))$, the distance away from origin equals $\sqrt{c_{1}^2+ \cdots + c_{n-1}^2+nt^2}$.

\medskip\noindent 

For the proof of (2) and (3), since the factors $(x_{i}^2-x_{j}^2)$ and $(x_{i}^2-t^2)$ remain unchanged along every fixed secondary slicing set $H_{c}$, then we only need to consider:

\medskip\noindent 

(i): For $C(n,n,n-1)(n\geq 2)$ and ${\bf C}(2n,2n-2)(n \geq 2)$, 
whether there exists radial integrals: $f(x)=g(|x|)$ such that

\begin{equation}\label{comf1}
 \frac{g\left(\sqrt{c_{1}^2+ \cdots + c_{n-1}^2+nt^2}\right)}{\prod\limits_{i=1}^{n-1}\sqrt{c_{i}^2+t^2}\sqrt{1+t^2\sum\limits_{i=1}^{n-1}\frac{1}{c_{i}^2+t^2}}} \geq \frac{g\left(\sqrt{c_{1}^2+ \cdots + c_{n-1}^2}\right)}{\prod\limits_{i=1}^{n-1}c_{i}},
\end{equation}
for all the value of $t$ and all the  secondary slicings $H_{c}$. 

\medskip\noindent 

(ii): For $C(n,n+1,n-1)(n\geq 2)$, whether there exists radial integrals: $f(x)=g(|x|)$ such that

\begin{equation}\label{comf2}
 \frac{g\left(\sqrt{c_{1}^2+ \cdots + c_{n-1}^2+nt^2}\right)}{\sqrt{1+t^2\sum\limits_{i=1}^{n-1}\frac{1}{c_{i}^2+t^2}}} \geq g\left(\sqrt{c_{1}^2+ \cdots + c_{n-1}^2}\right),
\end{equation}
for all the value of $t$ and all the  secondary slicings $H_{c}$. 

\medskip\noindent 

Compare \eqref{comf1} and \eqref{comf2}, easy to see: \eqref{comf1} has solutions for $g$ implies \eqref{comf2} has solutions for $g$; \eqref{comf2} doesn't have any solution for $g$ implies \eqref{comf1} has no solution for $g$.

\medskip\noindent 

Denote the symmetric polynomials in the $n-1$ variables: $c_{1}^2, \cdots, c_{n-1}^2$ by 
$$
\sigma_0=1,\sigma_1=\sum_{i=1}^{n-1}c_i^2,\ldots, \sigma_{n-2}=\sum_{j=1}^{n-1}c_1^2 \ldots \hat{c_j^2} \ldots c_{n-1}^2, \sigma_{n-1}=c_1^{2}\ldots c_{n-1}^2,
$$
since $f(x)=g(|x|)$ is non-negative,  \eqref{comf1} is equivalent to 

\begin{equation}\label{comf3}
g^2(\sqrt{\sigma_1+nt^2})\sigma_{n-1} \geq g^2(\sqrt{\sigma_1}) \left(\sum_{m=0}^{n-1}(m+1)t^{2m}\sigma_{n-1-m}\right)
\end{equation}

\medskip\noindent

\eqref{comf2} can be written as:

\begin{equation}\label{comf4}
g^2(\sqrt{\sigma_1+nt^2}) \geq g^2(\sqrt{\sigma_1})  \left(1+t^2\sum_{i=1}^{n-1}\frac{1}{c_{i}^2+t^2}\right) 
\end{equation}

\medskip\noindent

(i): $n=2$, \eqref{comf3} is equivalent to 
$$
g^2(\sqrt{\sigma_1+2t^2})\sigma_{1} \geq g^2(\sqrt{\sigma_1})(\sigma_1+2t^2),
$$
it is satisfied iff $\frac{g(r)}{r}$ is non-decreasing, some reasonable chosen integrals are: $f(x)=g(|x|)=|x|^{\alpha}(\alpha \geq 1)$,that prove (2).

\medskip\noindent 

(ii): $n\geq 3$, if \eqref{comf4} is satisfied for any secondary slicing sets $H_{c}$ and any points belong to $H_{c}$, then by firstly fixing $\sigma_1$ and $t^2<<\sigma_1$, and let $c_{n-1}<<t^2$ tends to $0$ (i.e. $H_c$ is very close to the focal point sets of the cone. Note that the focal distance is $|t|=x_{n-1}$, see \cite[Corollary 4.8]{CJX24}, etc), we can find that 
$$
g^2(\sqrt{\sigma_1+nt^2}) \geq 2 g^2(\sqrt{\sigma_1}), 
$$
this leads that the derivative of $g^2$ is infinity, a contradiction. And, we have proved (3). $ \Box $

\medskip\noindent 
\medskip\noindent

\textbf{4.2.3. Proof of Theorem 1.12:} Let $S$ be a rectifiable current having the same boundary with $\mathcal{C}_1$ or a rectifiable current reduced mod 2 having the same mod-2 boundary with $\mathcal{C}_1$. Denote the primary slicing sets of $\mathcal{C}_{1}$ by $T_{a}(a\in \mathcal{A})$, their union is denoted by $T$. The associated slices of $\mathcal{C}_{1}$ and $S$ are denoted by $C_{a}$ and $S_{a}$: $C_{a}=T_{a}\cap \mathcal{C}_{1}$, $S_{a}=T_{a}\cap S$. Recall that $T_{a}(a\in \mathcal{A})$ are disjoint normal space slicings.

We first prove that the following inequality is satisfied: 
$$
\int_{\mathcal{C}_1\cap H}fw_{1}\leq \int_{S\cap H }fw_{1},
$$
then by the isometric action, we see the inequality holds for every other primary slicing set $C_{a}$ and $S_{a}$ since $S$ is chosen arbitrarily.

\medskip\noindent

We can only consider the level sets of $h$ for that $H_{c}$ intersects the cones within the open unit ball, let $U$ be their union. $U$ is open in $H$ and covers $\mathcal{C}_1\cap H$ almost everywhere, by the area-coarea formula, note that every slicing set $H_{c}\subset U$ intersects $\mathcal{C}_1$ orthogonally, by area-coarea formula, note the weighting function $w_2$ are the reciprocal of the Jacobian of $h$, hence we have

$$
\int_{\mathcal{C}_1\cap H}fw_{1}=\int_{\mathcal{C}_1\cap U}fw_{1}= \int_{c \in \textbf{R}^{r}}fw_{1}w_{2}(H_{c}\cap \mathcal{C}_1) 
$$

Similar to the proof of Theorem 1.10, $S$ intersects almost every slicing set $H_{c}\subset U$. Now, combined with Proposition 4.1, it deduces that 

$$
\int_{c \in \textbf{R}^{r}}fw_{1}w_{2}(H_{c}\cap \mathcal{C}_1)\leq 
\int_{c \in \textbf{R}^{r}} \sum_{x\in H_{c}\cap S} fw_{1}w_{2}(x) \leq \int_{S\cap H}fw_{1},
$$
the first inequality follows from Proposition 4.1, the last inequality is again by the area-coarea formula, the sum indices arise from the fact that every $H_{c}$ intersects $S$ at least one point, maybe more! 

\medskip\noindent 

Then 
$$
\int_{\mathcal{C}_1\cap H}fw_{1} \leq \int_{S\cap H}fw_{1}.
$$

Having this inequality, let $S$ be a rectifiable current with the same boundary as $\mathcal{C}_1$, note that $T$ covers $\mathcal{C}_1$ a.e., then  
$$
\mathcal{E}_f(\mathcal{C}_1)=\mathcal{E}_f(\mathcal{C}_1\cap T)=\int_{a\in \mathcal{A}}\int_{C_{a}}fw_{1}\leq \int_{a\in \mathcal{A}}\int_{S_{a}}fw_{1} \leq \mathcal{E}_f(S \cap T) \leq \mathcal{E}_f(S),
$$ 
the second equality and the second inequality are both by the area-coarea formula. The final inequality follows the fact that $f$ is non-negative, this completes the proof. $\Box$

\medskip\noindent 
\medskip\noindent

\textbf{Acknowledgements} This work is supported by the NSFC (No. 12301068, No. 11871445), the project of Stable Support for Youth Team in Basic Research Field, CAS (YSBR-001), the China Postdoctoral Science Foundation (No. 2023M733401) and the Fundamental Research Funds for the Central Universities.

		\medskip\noindent
		\medskip\noindent

		\begin{flushleft}
			\medskip\noindent
			\begin{tabbing}
				XXXXXXXXXXXXXXXXXXXXXXXXXX*\=\kill
				Hongbin Cui\\
				School of Mathematical Sciences, University of Science and Technology of China\\
				Wu Wen-Tsun Key Laboratory of Mathematics, USTC, Chinese Academy of Sciences\\
				96 Jinzhai Road, Hefei, 230026, Anhui Province, China\\
				
				E-mail: cuihongbin@ustc.edu.cn
				
			\end{tabbing}
		\end{flushleft}

  \begin{flushleft}
			\medskip\noindent
			\begin{tabbing}
				XXXXXXXXXXXXXXXXXXXXXXXXXX*\=\kill
				Xiaowei Xu\\
				School of Mathematical Sciences, University of Science and Technology of China\\
				Wu Wen-Tsun Key Laboratory of Mathematics, USTC, Chinese Academy of Sciences\\
				96 Jinzhai Road, Hefei, 230026, Anhui Province, China\\
				
				E-mail: xwxu09@ustc.edu.cn
				
			\end{tabbing}
		\end{flushleft}
		
  \end{document}